\documentclass[12pt,a4paper]{article}%
\usepackage{amsmath}
\usepackage{graphicx}%
\usepackage{amsfonts}%
\usepackage{amssymb}

\begin{document}

\title{Morse functions and cohomology of homogeneous spaces}
\author{Haibao Duan\thanks{Supported by Polish KBN grant No.2 P03A 024 23.}\\Institute of Mathematics, \\Chinese Academy of Sciences, Beijing 100080\\dhb@math.ac.cn}
\date{\ }
\maketitle

This article arose from a series of three lectures given at the Banach Center,
Warsaw, during period of 24 March to 13 April, 2003.

\bigskip

Morse functions are useful tool in revealing the geometric formation of its
domain manifolds $M$. They define the handle decompositions of $M$ from which
the additive homologies $H_{\ast}(M)$ may be constructed. In these lectures
two further questions were emphasized.

\begin{quote}
(1) How to find a Morse function on a given manifold?

(2) From Morse functions can one derive the multiplicative cohomology rather
than the additive homology?
\end{quote}

\noindent It is not our intention here to make detailed studies of these
question. Instead, we will illustrate by examples solutions to them for some
classical manifolds as homogeneous spaces.

\bigskip

I am very grateful to Piotr Pragacz for the opportunity to speak of the wonder
that I have experienced with Morse functions, and for his hospitality during
my stay in Warsaw. Thanks are also due to Dr. Marek Szyjewski for taking the
lecture notes from which the present article was initiated, and to Dr. M.
Borodzik for many improvements on the earlier version of the note.

\section{Computing homology: a classical method}

There are many ways to introduce Morse Theory. However, I would like to
present it in the effective computation of homology (cohomology) of manifolds.

Homology (cohomology) theory is a bridge between geometry and algebra in the
sense that it assigns to a manifold $M$ a graded abelian group $H_{\ast}(M)$
(graded ring $H^{\ast}(M)$), assigns to a map $f:M\rightarrow N$ between
manifolds the induced homomorphism

\begin{center}
$f_{\ast}:H_{\ast}(M)\rightarrow H_{\ast}(N)$ (resp. $f^{\ast}:H^{\ast
}(N)\rightarrow H^{\ast}(M)$).
\end{center}

\noindent During the past century this idea has been widely applied to
translate geometric problems concerning manifolds and maps between them to
problems about groups (or rings) and homomorphisms, so that by solving the
latter in the well-developed framework of algebra, one obtains solutions to
the problems initiated from geometry.

The first problem one encounters when working with homology theory is the
following one.

\textbf{Problem 1.} \textsl{Given a manifold }$M$\textsl{, compute }$H_{\ast
}(M)$\textsl{ (as a graded abelian group) and }$H^{\ast}(M)$\textsl{ (as a
graded ring).}

We begin by recalling a classical method to approach the additive homology of manifolds.

\begin{center}
\textbf{1--1. Homology of a cell complex}
\end{center}

The simplest geometric object in dimension $n$, $n\geq0$, is the unit ball
$D^{n}=\{x\in\mathbb{R}^{n}\mid\left\|  x\right\|  ^{2}\leq1\}$ in the
Euclidean $n$-space $\mathbb{R}^{n}=\{x=(x_{1},\cdots,x_{n})\mid x_{i}%
\in\mathbb{R}\}$, which will be called the \textsl{n-dimensional disk} (or
\textsl{cell})\textsl{ }. Its boundary presents us the simplest closed $(n-1)$
dimensional manifold, the $(n-1)$ \textsl{sphere}: $S^{n-1}=\partial
D^{n}=\{x\in\mathbb{R}^{n}\mid\left\|  x\right\|  ^{2}=1\}$.

\bigskip

\ \ \ \
{\includegraphics[
natheight=2.268400in,
natwidth=6.002700in,
height=1.8913in,
width=4.9165in
]%
{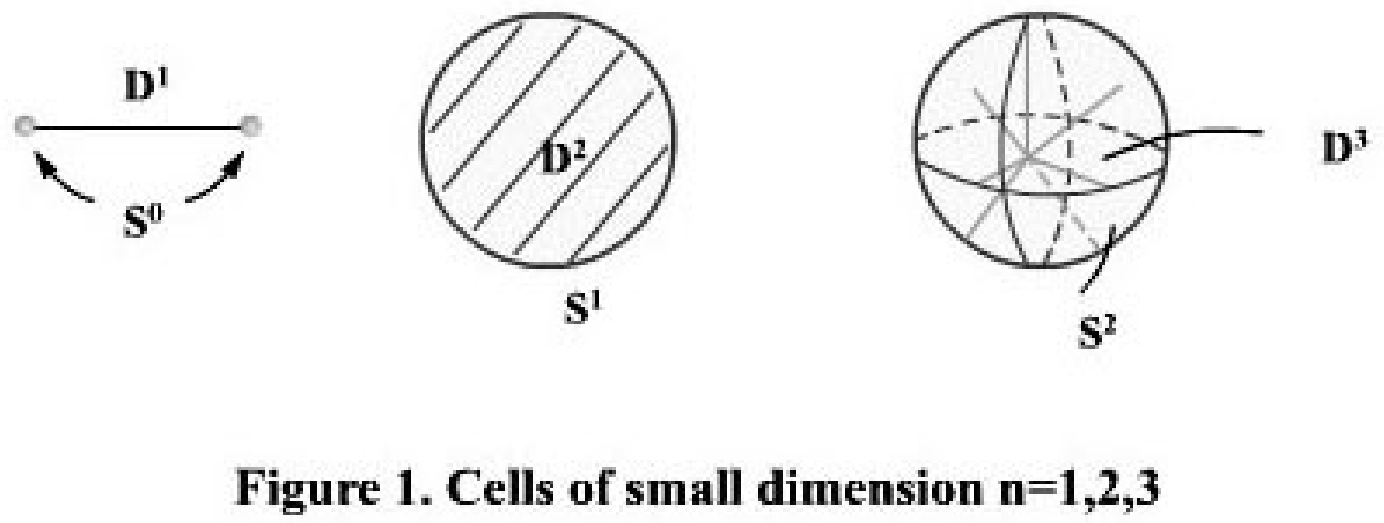}%
}%

\bigskip

Let $f:S^{r-1}\rightarrow X$ be a continuous map from $S^{r-1}$ to a
topological space $X$. Using $f$ we define

\begin{quote}
(1) an adjunction space
\end{quote}

\begin{center}
$X_{f}=X\cup_{f}D^{r}=X\sqcup D^{r}/y\in S^{r-1}\sim f(y)\in X$,
\end{center}

\begin{quote}
called \textsl{the space obtained from }$X$\textsl{ by attaching an n-cell
using }$f$.\ \

\
{\includegraphics[
natheight=2.415400in,
natwidth=6.288100in,
height=1.5298in,
width=4.3811in
]%
{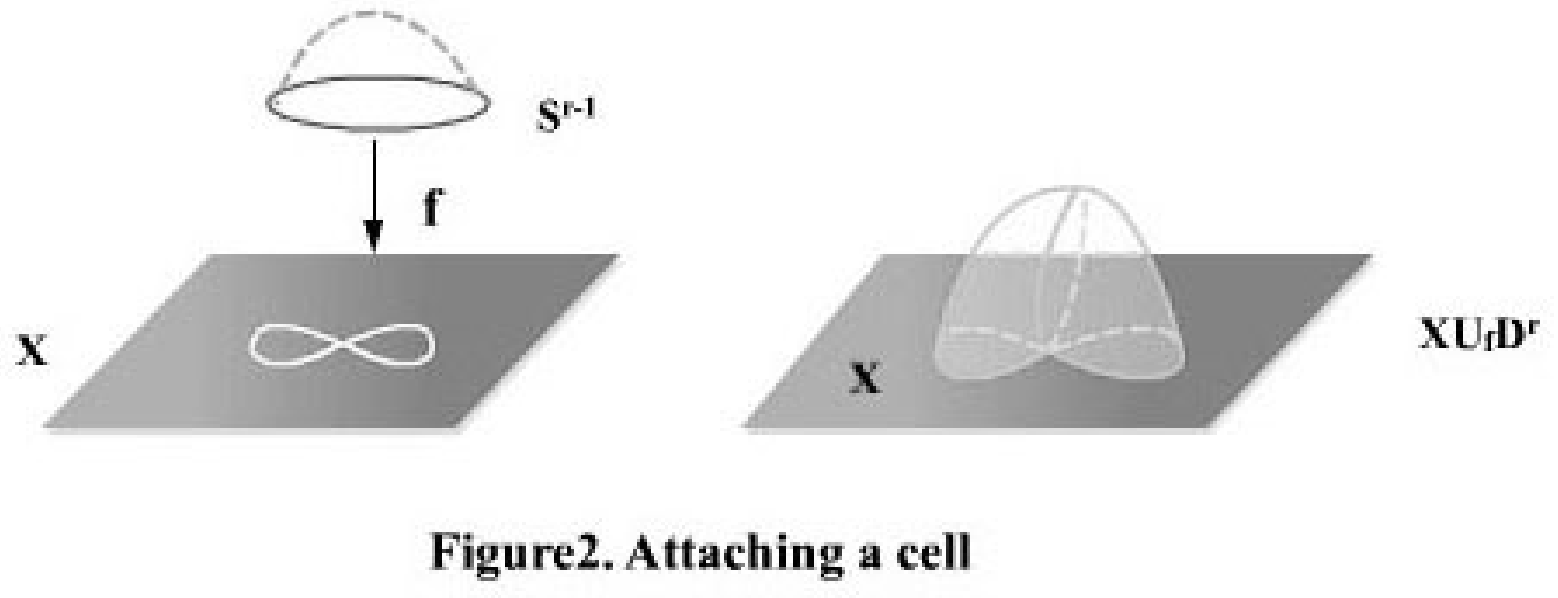}%
}%

(2) a homology class $f_{\ast}[S^{r-1}]\in H_{r-1}(X;\mathbb{Z})$ which
generates a cyclic subgroup of $H_{r-1}(X;\mathbb{Z})$: $a_{f}=<f_{\ast
}[S^{r-1}]>\subset H_{r-1}(X;\mathbb{Z})$.
\end{quote}

\noindent We observe that the integral homology of the new space $X\cup
_{f}D^{r}$ can be computed in terms of $H_{\ast}(X;\mathbb{Z})$ and its
subgroup $a_{f}$.

\textbf{Theorem 1.}\textsl{ Let }$X_{f}=X\cup_{f}D^{r}$\textsl{. Then the
inclusion }$i:X\rightarrow X_{f}$

\bigskip\textsl{1) induces isomorphisms }$H_{k}(X;\mathbb{Z})\rightarrow
H_{k}(X_{f};\mathbb{Z})$\textsl{ for all }$k\neq r,r-1$\textsl{;}

\textsl{2) fits into the short exact sequences}

\begin{center}
$0\rightarrow a_{f}\rightarrow H_{r-1}(X;\mathbb{Z})\overset{i_{\ast}%
}{\rightarrow}H_{r-1}(X_{f};\mathbb{Z})\rightarrow0$

$0\rightarrow H_{r}(X;\mathbb{Z})\overset{i_{\ast}}{\rightarrow}H_{r}%
(X_{f};\mathbb{Z})\rightarrow\{%
\begin{array}
[c]{c}%
0\text{ if }\left|  a_{f}\right|  =\infty\qquad\\
\mathbb{Z}\rightarrow0\text{ if }\left|  a_{f}\right|  <\infty\text{.}%
\end{array}
$
\end{center}

\textbf{Proof.} Substituting in the homology exact sequence of the pair
$(X_{f},X)$

\begin{center}
$H_{k}(X_{f},X;\mathbb{Z})=\{%
\begin{array}
[c]{c}%
0\text{ if }k\neq r;\\
\mathbb{Z}\text{ if }k=r
\end{array}
$
\end{center}

\noindent(note that the boundary operator maps the generator of $H_{r}%
(X_{f},X;\mathbb{Z})=\mathbb{Z}$ to $f_{\ast}[S^{r-1}]$), one obtains (1) and
(2) of the Theorem.$\square$

\bigskip

\textbf{Definition 1.1.} \textsl{Let }$X$\textsl{ be a topological space. A
cell-decomposition of }$X$\textsl{ is a sequence of subspaces }$X_{0}\subset
X_{1}\subset\cdots\subset X_{m-1}\subset X_{m}=X$\textsl{ so that}

\textsl{a) }$X_{0}$\textsl{ consists of finite many points }$X_{0}%
=\{p_{1},\cdots,p_{l}\}$\textsl{; and}

\textsl{b) }$X_{k}=X_{k-1}\cup_{f_{i}}D^{r_{k}}$\textsl{, where }%
$f_{i}:\partial D^{r_{k}}=S^{r_{k}-1}\rightarrow X_{k-1}$\textsl{ is a
continuous map.}

\textsl{Moreover, }$X$\textsl{ is called a cell complex if a
cell-decomposition of }$X$\textsl{ exists.}

\bigskip

Two comments are ready for the notion of cell-complex $X$.

\begin{quote}
(1) It can be build up using the simplest geometric objects $D^{n}$,
$n=1,2,\cdots$ by repeated applying the same construction as ``attaching cell'';

(2) Its homology can be computed by repeated applications of the single
algorithm (i.e. Theorem 1).
\end{quote}

The concept of cell-complex was initiated by Ehresmann in 1933-1934. Suggested
by the classical work of H. Schubert in algebraic geometry in 1879 [Sch],
Ehresmann found a cell decomposition for the complex Grassmannian manifolds
from which the homology of these manifolds were computed [Eh]. The cells
involved are currently known as \textsl{Schubert cells }(\textsl{varieties}) [MS].

In 1944, Whitehead [Wh] described a cell decomposition for the real Stiefel
manifolds (including all real orthogonal groups) in order to compute the
homotopy groups of these manifolds, where the cells were called \textsl{the
normal cells} by Steenrod [St] or \textsl{Schubert cells} by Dieudonn\'{e} [D,
p.226]. In terms of this cell decompositions the homologies of these manifolds
were computed by C. Miller in 1951 [M]. We refer the reader to Steenrod [St]
for the corresponding computation for complex and quaternionic Stiefel manifolds.

Historically, finding a cell decomposition of a manifold was a classical
approach to computing its homology. It should be noted that it is generally a
difficult and tedious task to find (or to describe) a cell-decomposition for a
given manifold. We are looking for simpler methods.

\begin{center}
\textbf{1--2. Attaching handles (Construction in manifolds)}
\end{center}

``Attaching cells'' is a geometric procedure to construct topological spaces
by using the elementary geometric objects $D^{r}$, $r\geq0$. The corresponding
construction in manifolds are known as ``attaching handles''or more
intuitively, ``\textsl{attaching thickened cells}''.

Let $M$ be an $n$-manifold with boundary $N=\partial M$, and let
$f:S^{r-1}\rightarrow N$ be a smooth embedding of an $(r-1)$-sphere whose
tubular neighborhood in $N$ is trivial: $T(S^{r-1})=S^{r-1}\times D^{n-r}$. Of
course, as in the previous section, one may form a new topological space
$M_{f}=M\cup_{f}D^{r}$ by attaching an $r$-cell to $M$ by using $f$. However,
the space $M_{f}$ is in general not a manifold!

\ \ \ \
{\includegraphics[
natheight=3.709200in,
natwidth=6.915000in,
height=2.335in,
width=4.2471in
]%
{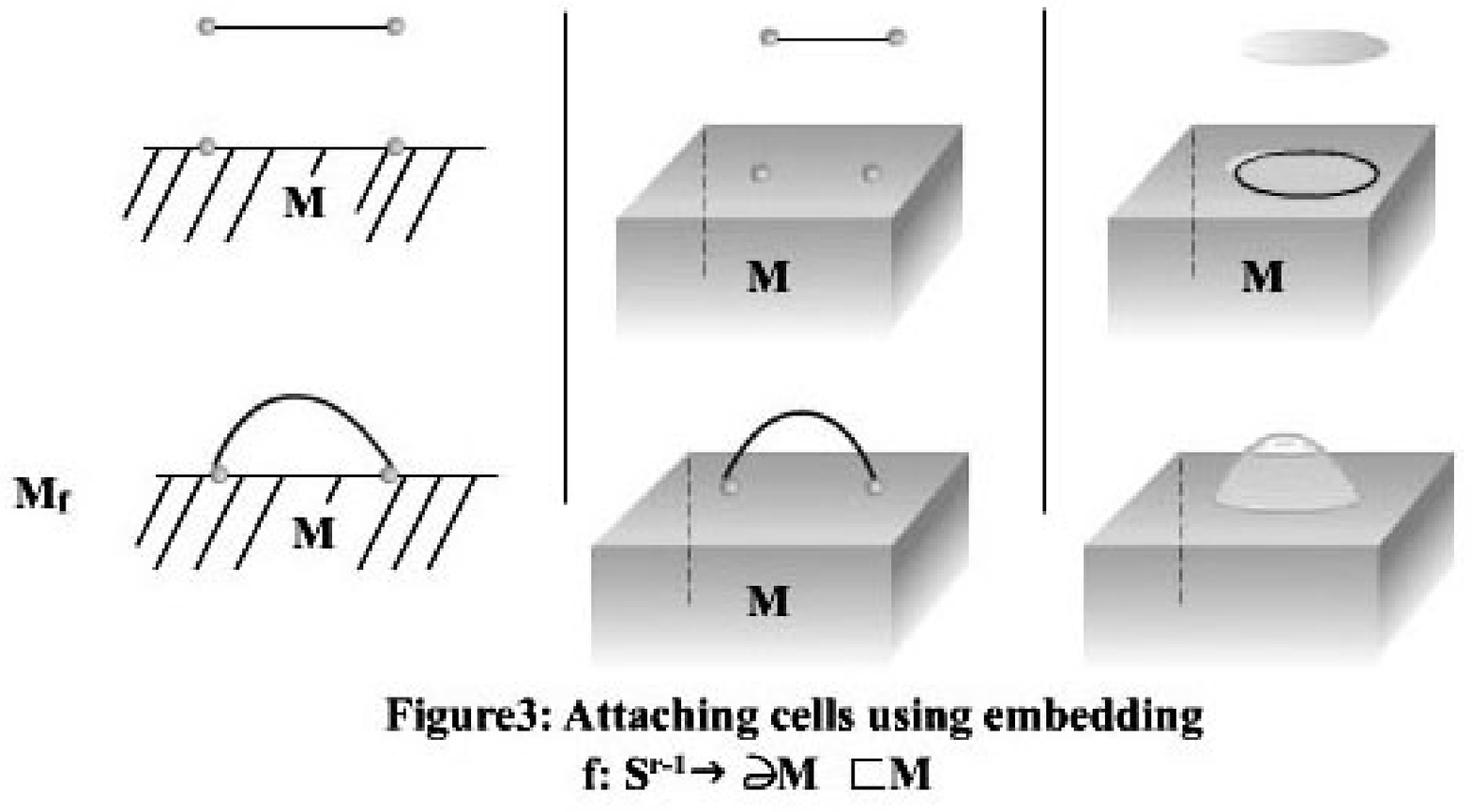}%
}%

Nevertheless, one may construct a new manifold $M^{^{\prime}}$ which contains
the space $M_{f}$ as a ``\textsl{strong deformation retract}'' by the
procedure below.

\textbf{Step 1.} To match the dimension of $M$, thicken the $r$-disc $D^{r}$
by taking product with $D^{n-r}$

\begin{center}
$D^{r}\times0\subset D^{r}\times D^{n-r}$ (a thickened $r$-disc)
\end{center}

\noindent and note that $\partial(D^{r}\times D^{n-r})=S^{r-1}\times
D^{n-r}\cup D^{r}\times S^{n-r-1}$.

\textbf{Step 2.} Choose a diffeomorphism

\begin{center}
$S^{r-1}\times D^{n-r}(\subset D^{r}\times D^{n-r})\overset{\varphi
}{\rightarrow}T(S^{r})\subset M$
\end{center}

\noindent that extends $f$ in the sense that $\varphi\mid S^{r-1}%
\times\{0\}=f$;

\textbf{Step 3.} Gluing $D^{r}\times D^{n-r}$ to $M$ by using $\varphi$ to
obtain $M^{^{\prime}}=M\cup_{\varphi}D^{r}\times D^{n-r}$.

\ \ \ \ \
{\includegraphics[
natheight=4.323200in,
natwidth=6.593300in,
height=2.5521in,
width=4.1701in
]%
{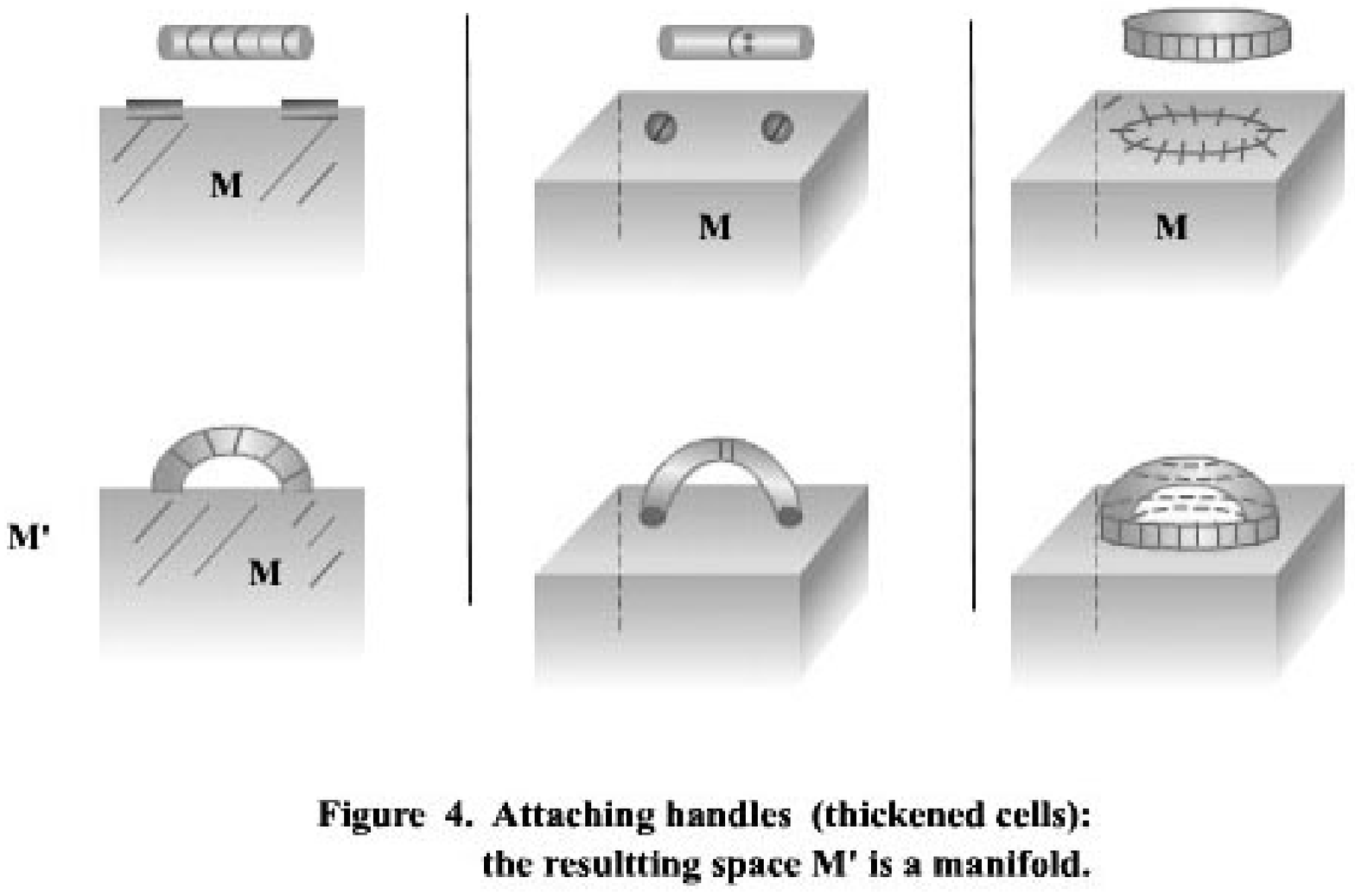}%
}%

\textbf{Step 4.} Smoothing the angles [M$_{3}$].

\bigskip

\textbf{Definition 1.2.} $M^{^{\prime}}$\textsl{ is called the manifold
obtained from }$M$\textsl{ by adding a thickened }$r-$\textsl{cell with core
}$M_{f}$\textsl{.}

\textbf{Remark.} The homotopy type (hence the homology) of $M^{^{\prime}}$
depends on the homotopy class $[f]\in\pi_{r-1}(M)$ of $f$.

The diffeomorphism type of $M^{^{\prime}}$ depends on the isotopy class of the
embedding $f$ (with trivial normal bundle), and a choice of $\varphi\in\pi
_{r}(SO(n-r))$.

\bigskip

Inside $M^{^{\prime}}=M\cup_{\varphi}D^{r}\times D^{n-r}$ one finds the
submanifold $M\subset M^{^{\prime}}$ as well as the subspace $M_{f}=M\cup
_{f}D^{r}\times\{0\}\subset M^{^{\prime}}=M\cup_{\varphi}D^{r}\times D^{n-r}$
in which the inclusion $j:M_{f}\rightarrow M^{^{\prime}}$ is a homotopy
equivalence. In particular, $j$ induces isomorphism in every dimension

\begin{center}
$H_{k}(M_{f},\mathbb{Z})\rightarrow H_{k}(M^{^{\prime}};\mathbb{Z})$, $k\geq0$.
\end{center}

\noindent Consequently, the integral cohomology of the new manifold
$M^{^{\prime}}$ can be expressed in terms of that of $M$ together with the
class $f_{\ast}[S^{r-1}]\in H_{r-1}(M;\mathbb{Z})$ by Theorem 1.

\textbf{Corollary.} \textsl{Let }$M^{^{\prime}}$\textsl{ be the manifold
obtained from }$M$\textsl{ by adding a thickened }$r-$\textsl{cell with core
}$M_{f}$\textsl{. Then the inclusion }$i:M\rightarrow M^{^{\prime}}$

\textsl{1) induces isomorphisms }$H_{k}(M;\mathbb{Z})\rightarrow
H_{k}(M^{^{\prime}};\mathbb{Z})$\textsl{ for all }$k\neq r,r-1$\textsl{;}

\textsl{2) fits into the short exact sequences}

\begin{center}
$0\rightarrow a_{f}\rightarrow H_{r-1}(M;\mathbb{Z})\rightarrow H_{r-1}%
(M^{^{\prime}};\mathbb{Z})\rightarrow0$

$0\rightarrow H_{r}(M;\mathbb{Z})\rightarrow H_{r}(M^{^{\prime}}%
;\mathbb{Z})\rightarrow\{%
\begin{array}
[c]{c}%
0\text{ if }\left|  a_{f}\right|  =\infty\qquad\\
\mathbb{Z}\rightarrow0\text{ if }\left|  a_{f}\right|  <\infty\text{.}%
\end{array}
$
\end{center}

\bigskip

\textbf{Definition 1.3.} \textsl{Let }$M$\textsl{ be a smooth closed }%
$n$\textsl{-manifold (with or without boundary). A handle decomposition of
}$M$\textsl{ is a filtration of submanifolds }$M_{1}\subset M_{2}\subset
\cdots\subset M_{m-1}\subset M_{m}=M$\textsl{ so that}

\textsl{(1) }$M_{1}=D^{n}$\textsl{;}

\textsl{(2) }$M_{k+1}$\textsl{ is a manifold obtained from }$M_{k}$\textsl{ by
attaching a thickened }$r_{k}$\textsl{-cell, }$r_{k}\leq n$\textsl{.}

\bigskip

If $M$ is endowed with a handle decomposition, its homology can be computed by
repeated applications of the corollary

\begin{center}
$H_{\ast}(M_{1})\mapsto H_{\ast}(M_{2})\mapsto\cdots\mapsto H_{\ast}(M)$.
\end{center}

\noindent Now, Problem 1 can be stated in geometric terms.

\textbf{Problem 2.} \textsl{Let }$M$\textsl{ be a smooth manifold.}

\textsl{(1) Does }$M$\textsl{ admits a handle decomposition?}

\textsl{(2) If yes, find one.}

\section{Elements of Morse Theory}

Using Morse function we prove, in this section, the following result which
answers (1) of Problem 2 affirmatively.

\textbf{Theorem 2. }\textsl{Any closed smooth manifold admits a handle decomposition.}

\begin{center}
\textbf{2--1. Study manifolds by using functions: the idea}
\end{center}

Let $M$ be a smooth closed manifold of dimension $n$ and let $f:M\rightarrow
\mathbb{R}$ be a non-constant smooth function on $M$. Put

\begin{center}
$a=\min\{f(x)\mid x\in M\}$, $b=\max\{f(x)\mid x\in M\}$.
\end{center}

\noindent Then $f$ is actually a map onto the interval $[a,b]$.

Intuitively, $f$ assigns to each point $x\in M$ a \textsl{height}
$f(x)\in\lbrack a,b]$. For a $c\in(a,b)$, those points on $M$ with the same
height $c$ (i.e. $L_{c}=f^{-1}(c)$) form \textsl{the level surface} of $f$ at
level $c$. It cuts the whole manifold into two parts $M=M_{c}^{-}\cup
M_{c}^{+}$ with

\begin{center}
$M_{c}^{-}=\{x\in M\mid f(x)\leq c\}$ (the part below $L_{c}$)

$M_{c}^{+}=\{x\in M\mid f(x)\geq c\}$ (the part above $L_{c}$)
\end{center}

\noindent and with $L_{c}=M_{c}^{-}\cap M_{c}^{+}$.

\ \ \
{\includegraphics[
natheight=2.977600in,
natwidth=4.440800in,
height=2.188in,
width=3.2223in
]%
{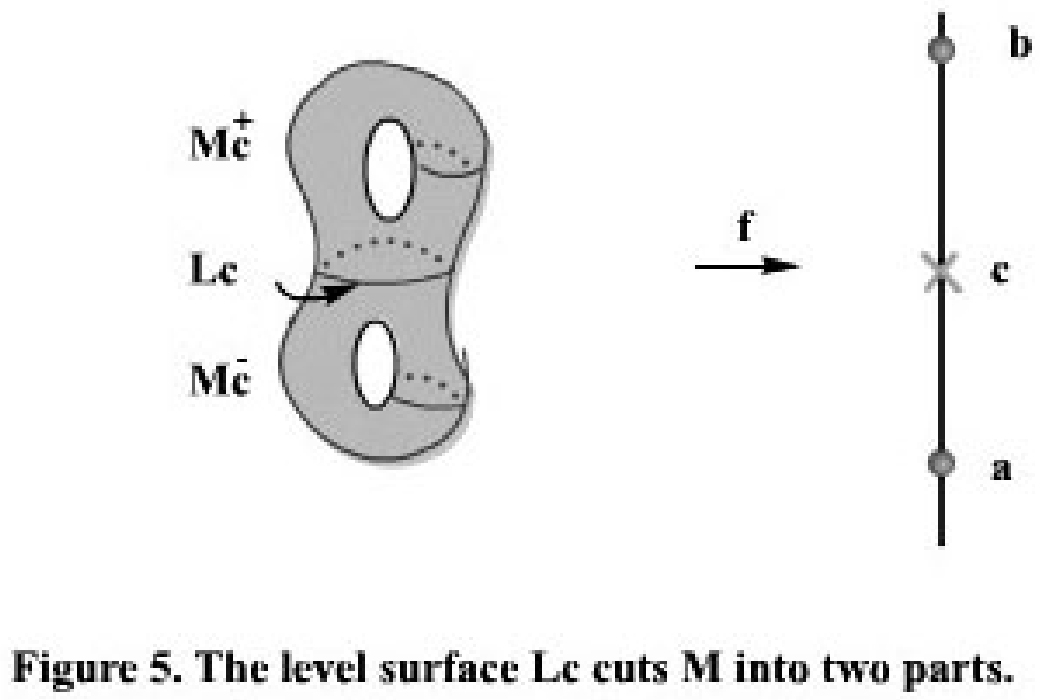}%
}%

In general, given a sequence of real numbers $a=c_{1}<\cdots<c_{m}=b$ , the
$m-2$ level surfaces $L_{c_{i}}$, $2\leq i\leq m-1$, defines a filtration on $M$

\begin{center}
$M_{1}\subset M_{2}\subset\cdots\subset M_{m-1}\subset M_{m}=M$,
\end{center}

\noindent with $M_{i}=M_{c_{i}}^{-}$.

\bigskip

Our aim is to understand the geometric construction of $M$ (rather than the
functions on $M$). Naturally, one expects to find a good function $f$ as well
as suitable reals $a=c_{1}<c_{2}<\cdots<c_{m}=b$ so that

\begin{quote}
(1) each $M_{i}$ is a smooth manifold with boundary $L_{c_{i}}$;

(2) the change in topology between each adjoining pair $M_{k}\subset M_{k+1}$
is as simple as possible.
\end{quote}

\noindent If this can be done, we may arrive at a global picture of the
construction of $M$.

Among all smooth functions on $M$, Morse functions are the most suitable for
this purpose.

\begin{center}
\textbf{2--2. Morse functions}
\end{center}

Let $f:M\rightarrow\mathbb{R}$ be a smooth function on a $n$-dimensional
manifold $M$ and let $p\in M$ be a point. In a local coordinates
$(x_{1},\cdots,x_{n})$ centered at $p$ (i.e. a Euclidean neighborhood around
$p$) the Taylor expansion of $f$ near $p$ reads

\begin{center}
$f(x_{1},\cdots,x_{n})=a+\underset{1\leq i\leq n}{\Sigma}b_{i}x_{i}%
+\underset{1\leq i,j\leq n}{\Sigma}c_{ij}x_{i}x_{j}+o(\parallel x\parallel
^{3})$,
\end{center}

\noindent in which

\begin{center}
$a=f(0)$;$\qquad b_{i}=\frac{\partial f}{\partial x_{i}}(0)$, $1\leq i\leq n$; and

$c_{ij}=\frac{1}{2}\frac{\partial^{2}f}{\partial x_{j}\partial x_{i}}(0)$,
$1\leq i,j\leq n$.
\end{center}

\ \
{\includegraphics[
natheight=2.296100in,
natwidth=5.055700in,
height=1.7945in,
width=4.4019in
]%
{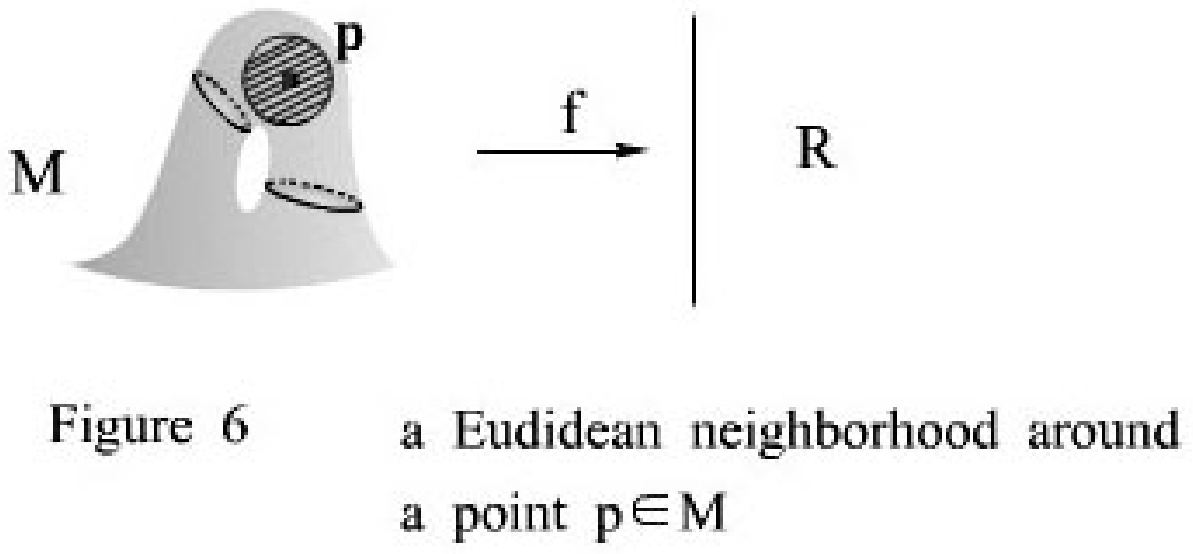}%
}%

\noindent Let $T_{p}M$ be the tangent space of $M$ at $p$. The $n\times n$
symmetric matrix,

\begin{center}
$H_{p}(f)=(c_{ij}):T_{p}M\times T_{p}M\rightarrow\mathbb{R}$ (resp.
$T_{p}M\rightarrow T_{p}M$)
\end{center}

\noindent called the \textsl{Hessian form} (resp. Hessian operator) of $f$ at
$p$, can be brought into diagonal form by changing the linear basis
$\{\frac{\partial}{\partial x_{1}},\cdots,\frac{\partial}{\partial x_{n}}\}$
of $T_{p}M$

\begin{center}
$H_{p}(f)=(c_{ij})\sim0_{s}\oplus(-I_{r})\oplus(I_{t})$, $s+r+t=n$.
\end{center}

\bigskip

\textbf{Definition 2.1.} $p\in M$\textsl{ is called a critical point of }%
$f$\textsl{ if in a local coordinates at }$p$\textsl{, }$b_{i}=0$\textsl{ for
all }$1\leq i\leq n$\textsl{. Write }$\Sigma_{f}$\textsl{ for the set of all
critical points of }$f$\textsl{.}

\textsl{A critical point }$p\in\Sigma_{f}$\textsl{ is called non-degenerate if
the form }$H_{p}(f)$\textsl{ is non-degenerate. In this case the number }%
$r$\textsl{ is called the index of }$p$\textsl{ (as a non-degenerate critical
point of }$f$\textsl{), and will be denoted by }$r=Ind(p)$\textsl{.}

$f$\textsl{ is said to be a Morse function on }$M$\textsl{ if all its critical
points are non-degenerate.}

\bigskip

The three items ``\textsl{critical point}'', ``\textsl{non-degenerate critical
point}'' as well as the ``\textsl{index}'' of a nondegenerate critical point
specified in the above are clearly independent of the choice of local
coordinates centered at $p$. Two useful properties of a Morse function are
given in the next two lemmas.

\textbf{Lemma 2.1.} \textsl{If }$M$\textsl{ is closed and if }$f$\textsl{ is a
Morse function on }$M$\textsl{, then }$\Sigma_{f}$\textsl{ is a finite set.}

\textbf{Proof.} The set $\Sigma_{f}$ admits an intrinsic description without
referring to local coordinate systems.

The tangent map $Tf:TM\rightarrow\mathbb{R}$ of $f$ gives rise to a cross
section $\sigma_{f}:M\rightarrow T^{\ast}M$ \ for the cotangent bundle
$\pi:T^{\ast}M\rightarrow M$. Let $\sigma:M\rightarrow T^{\ast}M$ be the zero
section of $\pi$. Then $\Sigma_{f}=\sigma_{f}^{-1}[\sigma(M)]$. $f$ is a Morse
function is equivalent to the statement that the two embeddings $\sigma
_{f},\sigma:M\rightarrow T^{\ast}M$ have transverse intersection.$\square$

\bigskip

\textbf{Lemma 2.2 }(\textbf{Morse Lemma}, cf. [H; p.146]). \textsl{If }$p\in
M$\textsl{ is a non-degenerate critical point of }$f$\textsl{ with index }%
$r$\textsl{, there exist local coordinates }$(x_{1},\cdots,x_{n})$\textsl{
centered at }$p$\textsl{ so that}

\begin{center}
$f(x_{1},\cdots,x_{n})=f(0)-\underset{1\leq i\leq r}{\Sigma}x_{i}%
^{2}+\underset{r<i\leq n}{\Sigma}x_{i}^{2}$
\end{center}

\noindent\textsl{(i.e. the standard nondegenerate quadratic function of index
}$r$\textsl{)}.

\textbf{Proof.} By a linear coordinate change we may assume that

\begin{center}
$(\frac{\partial^{2}f}{\partial x_{j}\partial x_{i}}(0))=(-I_{r}%
)\oplus(I_{n-r})$.
\end{center}

\noindent Applying the fundamental Theorem of calculus twice yields the expansion

\begin{enumerate}
\item[(A)] $\qquad\qquad\qquad f(x_{1},\cdots,x_{n})=f(0)+\underset{1\leq
i,j\leq n}{\Sigma}\underset{}{x_{i}x_{j}}b_{ij}(x)$
\end{enumerate}

\noindent in which

\begin{center}
$b_{ij}(x)=\int_{0}^{1}\int_{0}^{1}\frac{\partial^{2}\overline{f}}{\partial
x_{j}\partial x_{i}}(stx_{1},\cdots,stx_{n})dtds$.
\end{center}

\noindent The family of matrix $B(x)=(b_{ij}(x))$, $x\in U$, may be considered
as a smooth map

\begin{center}
$B:U\rightarrow\mathbb{R}^{\frac{n(n+1)}{2}}$(=the vector space of all
$n\times n$ symmetric matrices).
\end{center}

\noindent with $B(0)=$ $(-I_{r})\oplus(I_{n-r})$, where $U\subset M$ is the
Euclidean neighborhood centered at $p$. It follows that

``\textsl{there is a smooth map }$P:$\textsl{ }$U\rightarrow GL(n)$\textsl{ so
that in some neighborhood }$V$\textsl{ of }$0\in U$\textsl{,}

\begin{center}
$B(x)=P(x)\{(-I_{r})\oplus(I_{n-r})\}P(x)^{\tau}$\textsl{ and }$P(0)=I_{n}$''.
\end{center}

\noindent With this we infer from (A) that, for $x=(x_{1},\cdots,x_{n})\in V$

\begin{center}
$f(x)=f(0)+xB(x)x^{\tau}=f(0)+xP(x)\{(-I_{r})\oplus(I_{n-r})\}P(x)^{\tau
}x^{\tau}$.
\end{center}

\noindent It implies that if one makes the coordinate change

\begin{center}
$(y_{1},\cdots,y_{n})=(x_{1},\cdots,x_{n})P(x)$
\end{center}

\noindent on a neighborhood of $0\in U$ then one gets

\begin{center}
$f(y_{1},\cdots,y_{n})=f(0)-\underset{1\leq i\leq r}{\Sigma}y_{i}%
^{2}+\underset{r<i\leq n}{\Sigma}y_{i}^{2}$.$\square$

\bigskip

\textbf{2--3. Geometry of gradient flow lines}
\end{center}

The first set of information we can derive directly from a Morse function
$f:M\rightarrow\mathbb{R}$ consists of

\begin{quote}
(1) the set $\Sigma_{f}$ of critical points of $f$;

(2) the index function $Ind:\Sigma_{f}\rightarrow\mathbb{Z}$.
\end{quote}

\bigskip

Equip $M$ with a Riemannian metric so that the gradient field of $f\qquad$

\begin{center}
$v=grad(f):M\rightarrow TM$,
\end{center}

\noindent is defined. One of the very first thing that one learns from the
theory of ordinary differential equations is that, for each $x\in M$, there
exists a unique smooth curve $\varphi_{x}:\mathbb{R}\rightarrow M$ subject to
the following constraints

\begin{quote}
(1) the initial condition: $\varphi_{x}(0)=x$;

(2) the ordinary differential equation:$\frac{d\varphi_{x}(t)}{dt}%
=v(\varphi_{x}(t))$;

(3) $\varphi_{x}$ varies smoothly with respect to $x\in M$ in the sense that

``the map $\varphi:M\times\mathbb{R}\rightarrow M$ by $(x,t)\rightarrow
\varphi_{x}(t)$ is smooth and, for every $t\in\mathbb{R}$, the restricted
function $\varphi:M\times\{t\}\rightarrow M$ is a diffeomorphism.''

\ \ \ \ \
{\includegraphics[
natheight=2.769100in,
natwidth=5.556400in,
height=1.9268in,
width=4.0309in
]%
{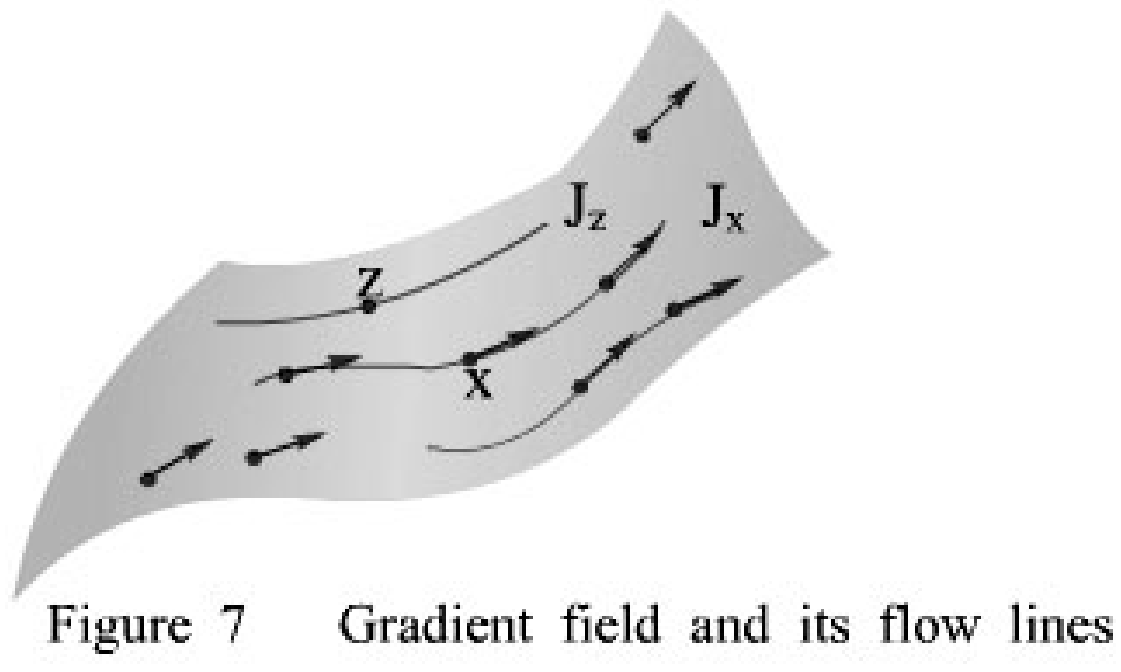}%
}%
\end{quote}

\bigskip

\textbf{Definition 2.2. }\textsl{For }$x\in M$\textsl{ \ let }$J_{x}%
=\operatorname{Im}$\textsl{ }$\varphi_{x}\subset M$\textsl{, and call it the
gradient flow line of }$f$\textsl{ through }$x$\textsl{.}

An alternative description for $J_{x}$ is the following. It is the image of
the parameterized curve $\varphi(t)$ in $M$ that satisfies

1) passing through $x$ at the time $t=0$;

2) at any point $y\in J_{x}$, the tangent vector $\frac{d\varphi}{dt}$ to
$J_{x}$ at $y$ agrees with the value of $v$ at $y$.

\bigskip

We build up the geometric picture of flow lines in the result below.

\textbf{Lemma 2.3 }(Geometry of gradient flow lines).

\textsl{(1) }$x\in\Sigma_{f}\Leftrightarrow J_{x}$\textsl{ consists of a point;}

\textsl{(2) }$\forall x,y\in M$\textsl{ we have either }$J_{x}=J_{y\text{ }}%
$\textsl{ or }$J_{x}\cap J_{y\text{ }}=\emptyset$\textsl{;}

\textsl{(3) if }$x\notin\Sigma_{f}$\textsl{, then }$J_{x}$\textsl{ meets level
surfaces of }$f$\textsl{ transversely; and }$f$\textsl{ is strictly increasing
along the directed curve }$J_{x}$\textsl{;}

\textsl{(4) if }$x\notin\Sigma_{f}$\textsl{, the two limits }$\underset
{t\rightarrow\pm\infty}{\lim}\varphi_{x}(t)$\textsl{ exist and belong to
}$\Sigma_{f}$\textsl{.}

\ \ \ \ \ \
{\includegraphics[
natheight=2.769100in,
natwidth=5.556400in,
height=2.0211in,
width=4.2808in
]%
{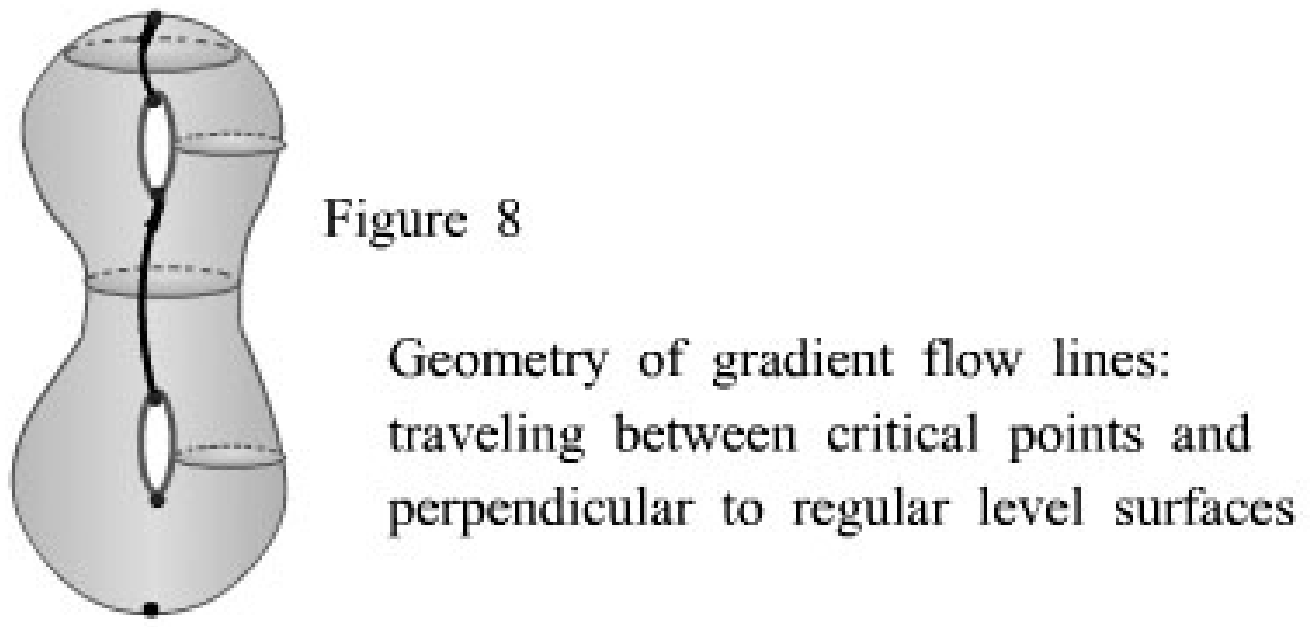}%
}%

\textbf{Proof}. (2) comes from the fact that $\varphi_{\varphi_{x}%
(t)}(s)=\varphi_{x}(t+s)$.

(3) is verified by

\begin{center}
$\frac{df\varphi_{x}(t)}{dt}=<$grad$f,\frac{d\varphi_{x}(t)}{dt}>=\mid
$grad$f\mid^{2}>0$.
\end{center}

Since the function $f\varphi_{x}(t)$ is bounded $a\leq f\varphi_{x}(t)\leq b$
and is monotone in $t$, the limits $\underset{t\rightarrow\pm\infty}{\lim}$
$f\varphi_{x}(t)$ exist. It follows from (3) that $\underset{t\rightarrow
\pm\infty}{\lim}\mid$grad$_{\varphi_{x}(t)}f\mid^{2}=0$. This shows
(4).$\square$

\bigskip

The most important notion subordinate to flow lines is:

\textbf{Definition 2.3.} \textsl{For a }$p\in\Sigma_{f}$\textsl{ we write}

\begin{center}
$S(p)=\underset{\underset{t\rightarrow+\infty}{\lim}\varphi_{x}(t)=p}{\cup
}J_{x}\cup\{p\}$\textsl{; }$\ T(p)=\underset{\underset{t\rightarrow-\infty
}{\lim}\varphi_{x}(t)=p}{\cup}J_{x}\cup\{p\}$\textsl{.}
\end{center}

\noindent\textsl{These will be called respectively the descending cell and the
ascending cell of }$f$\textsl{ at the critical point }$p$\textsl{.}

The term ``\textsl{cell}'' appearing in Definition 2.3 is justified by the
next result.

\textbf{Lemma 2.4.} \textsl{If }$p\in\Sigma_{f}$\textsl{ with Ind}%
$(p)=r$\textsl{, then }$(S(p),p)\cong(R^{r},0)$\textsl{, }$(T(p),p)\cong
(R^{n-r},0)$\textsl{, and both meet transversely at }$p$\textsl{.}

\textbf{Proof.} Let $(\mathbb{R}^{n},0)\subset(M,p)$ be an Euclidean
neighborhood centered at $p$ so that

\begin{center}
$f(x,y)=f(0)-\mid x\mid^{2}+\mid y\mid^{2}$(cf. Lemma 2.2),
\end{center}

\noindent where $(x,y)\in\mathbb{R}^{n}=\mathbb{R}^{r}\oplus\mathbb{R}^{n-r}$.
We first examine $S(p)\cap\mathbb{R}^{n}$ and $T(p)\cap\mathbb{R}^{n}$.

On $\mathbb{R}^{n}$ the gradient field of $f$ is easily seen to be
grad$f=(-2x,2y)$. The flow line $J_{x_{0}}$ through a point $x_{0}%
=(a,b)\in\mathbb{R}^{n}=\mathbb{R}^{r}\oplus\mathbb{R}^{n-r}$ is

\begin{center}
$\varphi_{x_{0}}(t)=(ae^{-2t},be^{2t}),t\in\mathbb{R}$.
\end{center}

\noindent Now one sees that

\begin{center}
$x_{0}\in S(p)\cap\mathbb{R}^{n}\Longleftrightarrow\underset{t\rightarrow
+\infty}{\lim}\varphi_{x_{0}}(t)=0$($p$)$\Longleftrightarrow b=0$;

$x_{0}\in T(p)\cap\mathbb{R}^{n}\Longleftrightarrow\underset{t\rightarrow
-\infty}{\lim}\varphi_{x_{0}}(t)=0$($p$)$\Longleftrightarrow a=0$.
\end{center}

\noindent It follows that

\begin{enumerate}
\item[(B)] $\qquad S(p)\cap\mathbb{R}^{n}=\mathbb{R}^{r}\oplus\{0\}\subset
\mathbb{R}^{n}$; $\quad T(p)\cap\mathbb{R}^{n}=\{0\}\oplus\mathbb{R}%
^{n-r}\subset\mathbb{R}^{n}$
\end{enumerate}

\noindent and both sets meet transversely at $0=p$.

\bigskip

\ \ \ \ \
{\includegraphics[
natheight=3.087400in,
natwidth=4.875000in,
height=2.226in,
width=3.314in
]%
{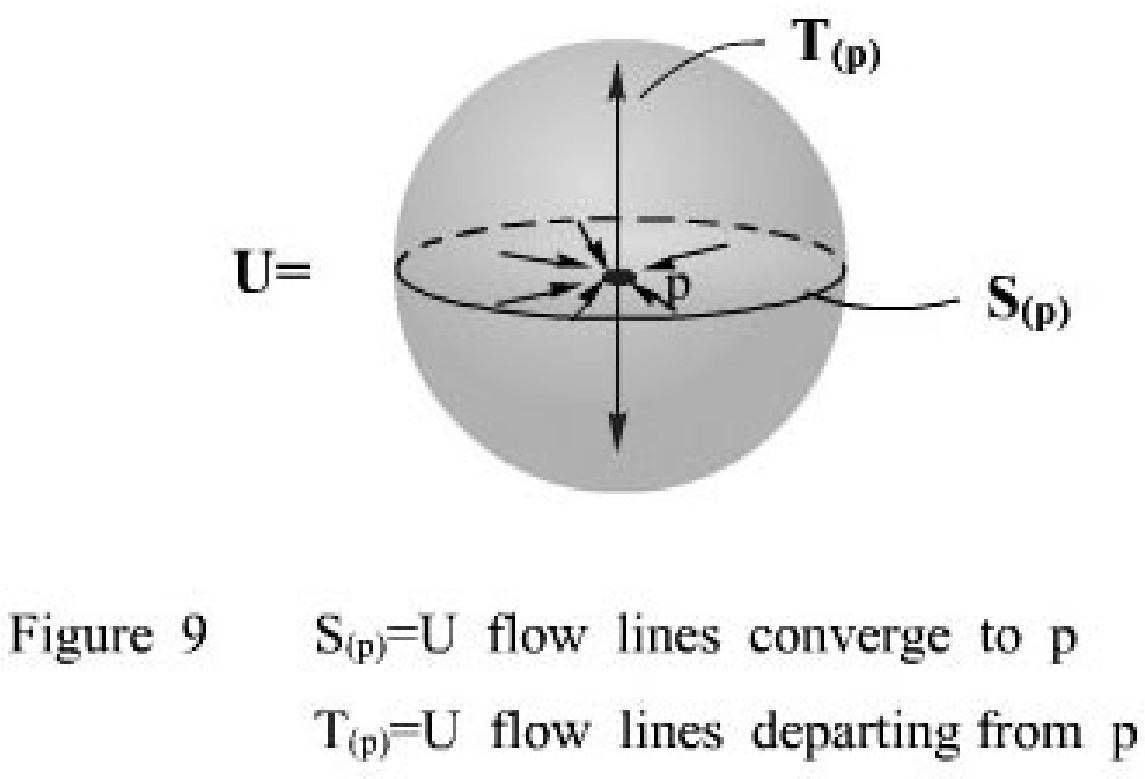}%
}%

Let $S^{n-1}$ be the unit sphere in $\mathbb{R}^{n}$ and put

\begin{center}
$S_{-}=S(p)\cap S^{n-1}$ (resp. $S_{+}=T(p)\cap S^{n-1}$).
\end{center}

\noindent Then (B) implies that $S_{-}\cong S^{r-1}$ (resp. $S_{+}\cong
S^{n-r-1}$). Furthermore, (2) of Lemma 2.3 implies that, for any $x\in S(p)$,
$J_{x}=J_{v}$ for some unique $v\in S_{-}$ because of $\varphi_{x}(t)\in
S(p)\cap\mathbb{R}^{n}$ for sufficient large $t$ with $\underset
{t\rightarrow+\infty}{\lim}\varphi_{x}(t)=p$ . Therefore

\begin{center}
$S(p)=\underset{v\in S_{-}}{\cup}J_{v}\cup\{p\}$ (resp. $T(p)=\underset{v\in
S_{+}}{\cup}J_{v}\cup\{p\}$).
\end{center}

\noindent That is, $S(p)$ (resp. $T(p)$) is an open cone over $S_{-}$ (resp.
$S_{+}$) with vertex $p$.$\square$

\bigskip

Summarizing, at a critical point $p\in\Sigma_{f}$,

\begin{quote}
(1) the flow lines that grow to $p$ (as $t\rightarrow\infty$) form an open
cell of dimension $Ind(p)=r$ centered at $p$ which lies below the critical
level $L_{f(p)}$;

(2) those flow lines that grow out of from $p$ (as $t\rightarrow\infty$) form
an open cell of dimension $Ind(p)=n-r$ centered at $p$ which lies above the
critical level $L_{f(p)}$.

\bigskip
\end{quote}

\begin{center}
\textbf{2--4. Handle decomposition of a manifold}
\end{center}

Our proof of Theorem 2 implies that the set of descending cells $\{S(p)\subset
M\mid p\in\Sigma_{f}\}$ of a Morse function on $M$ endows $M$ with the
structure of a cell complex.

\textbf{Proof of Theorem 2.} Let $f:M\rightarrow\lbrack a,b]$ be a Morse
function on a closed manifold $M$ with critical set $\Sigma_{f}$ and index
function $Ind:\Sigma_{f}\rightarrow\mathbb{Z}$. By Lemma 2.1 the set
$\Sigma_{f}$ is finite and we can assume that elements in $\Sigma_{f}$ are
ordered as $\{p_{1},\cdots,p_{m}\}$ by its values under $f$

\begin{center}
$a=f(p_{1})<f(p_{2})<\cdots<f(p_{m-1})<f(p_{m})=b$ [M$_{1}$, section 4].
\end{center}

\noindent Take a $c_{i}\in(f(p_{i}),f(p_{i+1}))$, $i\leq m-1$. Then $c_{i}$ is
a regular value of $f$. As a result $M_{i}=f^{-1}[a,c_{i}]\subset M$ is a
smooth submanifold with boundary $\partial M_{i}=L_{c_{i}}$. Moreover we get a
filtration on $M$ by submanifolds

\begin{center}
$M_{1}\subset M_{2}\subset\cdots\subset M_{m-1}\subset M_{m}=M$.
\end{center}

\noindent We establish theorem 2 by showing that

1) $M_{1}=D^{n}$;

2) For each $k$ there is an embedding $g:S^{r-1}\rightarrow\partial M_{k}$ so that

\begin{center}
$M_{k}\cup S(p_{k+1})=M_{k}\cup_{g}D^{r}$, $r=$Ind $(p_{k+1})$;
\end{center}

3) $M_{k+1}=M_{k}\cup D^{r}\times D^{n-r}$ with core $M_{k}\cup_{g}D^{r}$.

\bigskip

1) Let $\mathbb{R}^{n}$ be an Euclidean neighborhood around $p_{1}$ so that

\begin{center}
$f(x_{1},\cdots,x_{n})=a+\Sigma x_{i}^{2}$,
\end{center}

\noindent here we have made use of the fact $Ind(p_{1})=0$ (because $f$
attains its absolute minimal value $a$ at $p_{1}$) as well as Lemma 2.2. Since
$c_{1}=a+\varepsilon$ we have

\begin{center}
$f^{-1}[a,c_{1}]=\{x\in\mathbb{R}^{n}\mid\parallel x\parallel^{2}%
\leq\varepsilon\}\cong D^{n}$.
\end{center}

2) With the notation introduced in the proof of Lemma 2.4 we have

\begin{enumerate}
\item[(C)] $\qquad\qquad S(p_{k+1})=\underset{v\in S_{-}}{\cup}J_{v}%
\cup\{p_{k+1}\}$
\end{enumerate}

\noindent where $S_{-}\cong S^{r-1}$, $r=Ind(p_{k+1})$, and where $J_{v}$ is
the unique flow line $\varphi_{v}(t)$ with $\varphi_{v}(0)=v$ and with
$\underset{t\rightarrow+\infty}{\lim}\varphi_{v}(t)=$ $p_{k+1}$.

For a $v\in S_{-}$, $\underset{t\rightarrow-\infty}{\lim}\varphi_{v}%
(t)\in\{p_{1},\cdots,p_{k}\}\subset$ $Int(M_{k})$ by (4) and (3) of Lemma 2.3.
So $J_{v}$ must meet $\partial M_{k}$ at some unique point. The map
$g:S_{-}\rightarrow\partial M_{k}$ such that $g(v)=J_{v}\cap\partial M_{k}$ is
now well defined and must be an embedding by (2) of Lemma 2.3. We get

\begin{center}
$M_{k}\cup S(p_{k+1})=M_{k}\cup_{g}D^{r}$
\end{center}

\noindent form (C).

\bigskip

\ \ \ \ \
{\includegraphics[
natheight=3.882100in,
natwidth=5.556400in,
height=2.1473in,
width=3.2993in
]%
{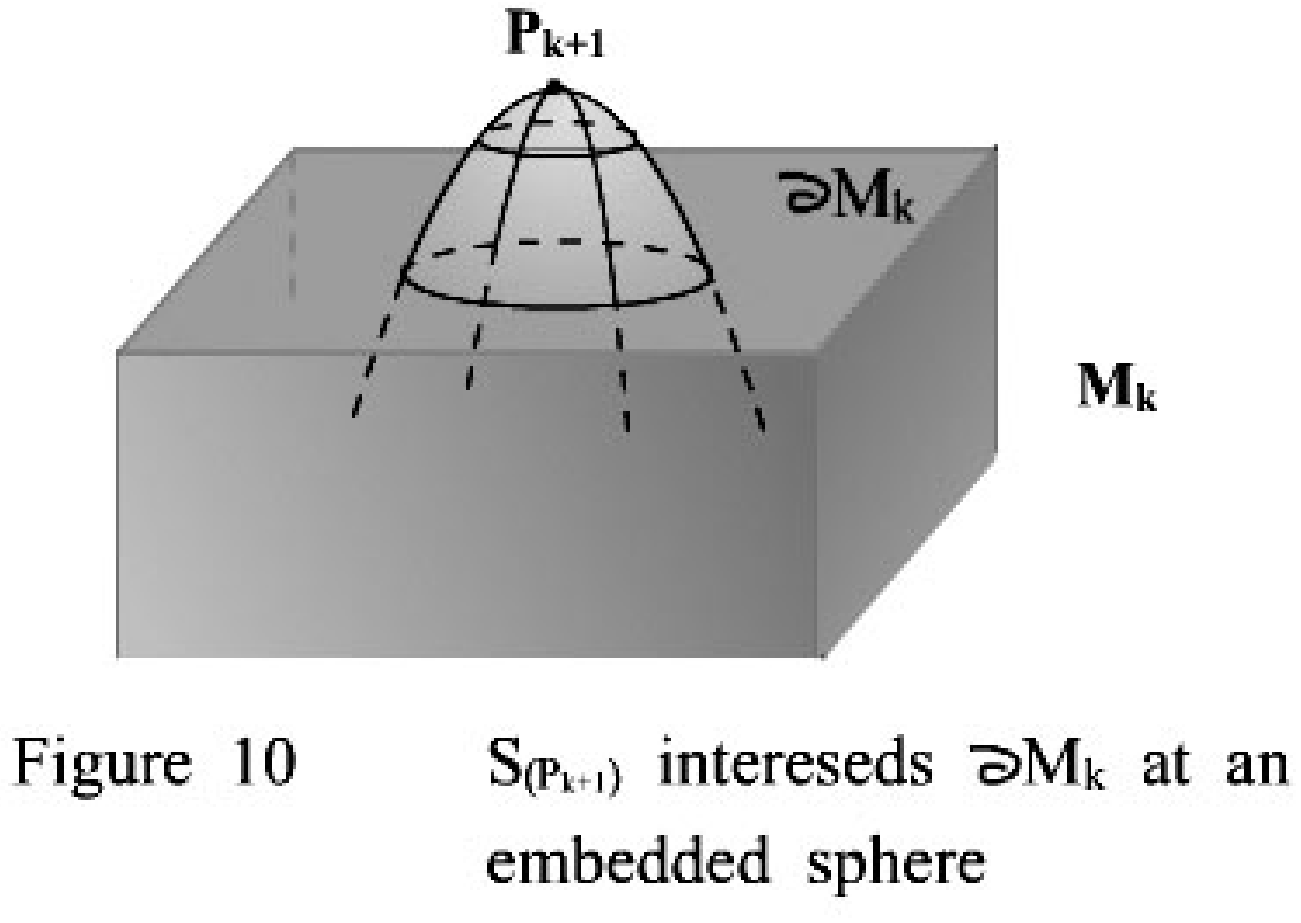}%
}%

3). In [M$_{1}$, p.33-34], Milnor demonstrated explicitly two deformation retractions

\begin{center}
$r:M_{k+1}\overset{R_{1}}{\rightarrow}M_{k}\cup D^{r}\times D^{n-r}%
\overset{R_{2}}{\rightarrow}M_{k}\cup S(p_{k+1})$
\end{center}

\noindent where $R_{1}$ does not change the diffeomorphism type of $M_{k+1}$
and where $D^{r}\times D^{n-r}$ is a thickening of the $r$-cell corresponding
to $S(p_{k+1})$.$\square$

\section{Morse functions via Euclidean geometry}

Our main subject is the effective computation of the additive homology or the
multiplicative cohomology of a given manifold $M$. Recall from section 1 that
if $M$ is endowed with a cell decomposition, the homology $H_{\ast}(M)$ can be
calculated by repeated application of Theorem 1. We have seen further in
section 2 that a Morse function $f$ on $M$ \ may define a cell-decomposition
on $M$ with each critical point of index $r$ corresponds to an $r$-cell in the
decomposition. The question that remains to us is

\begin{quote}
\textsl{How to find a Morse function on a given manifold?}
\end{quote}

\begin{center}
\textbf{3--1. Distance function on a Euclidean submanifold }
\end{center}

By a classical result of Whitney, every $n$-dimensional smooth manifold $M$
can be smoothly embedded into Euclidian space of some dimension less than
$2n+1$. Therefore, it suffices to assume that $M$ is a submanifold in an
Euclidean space $E$.

A point $a\in E$ gives rise to a function $f_{a}:M\rightarrow\mathbb{R}$ by
$f_{a}(x)=\parallel x-a\parallel^{2}$.

\bigskip

\ \ \
{\includegraphics[
natheight=3.019900in,
natwidth=5.528700in,
height=2.1266in,
width=3.7118in
]%
{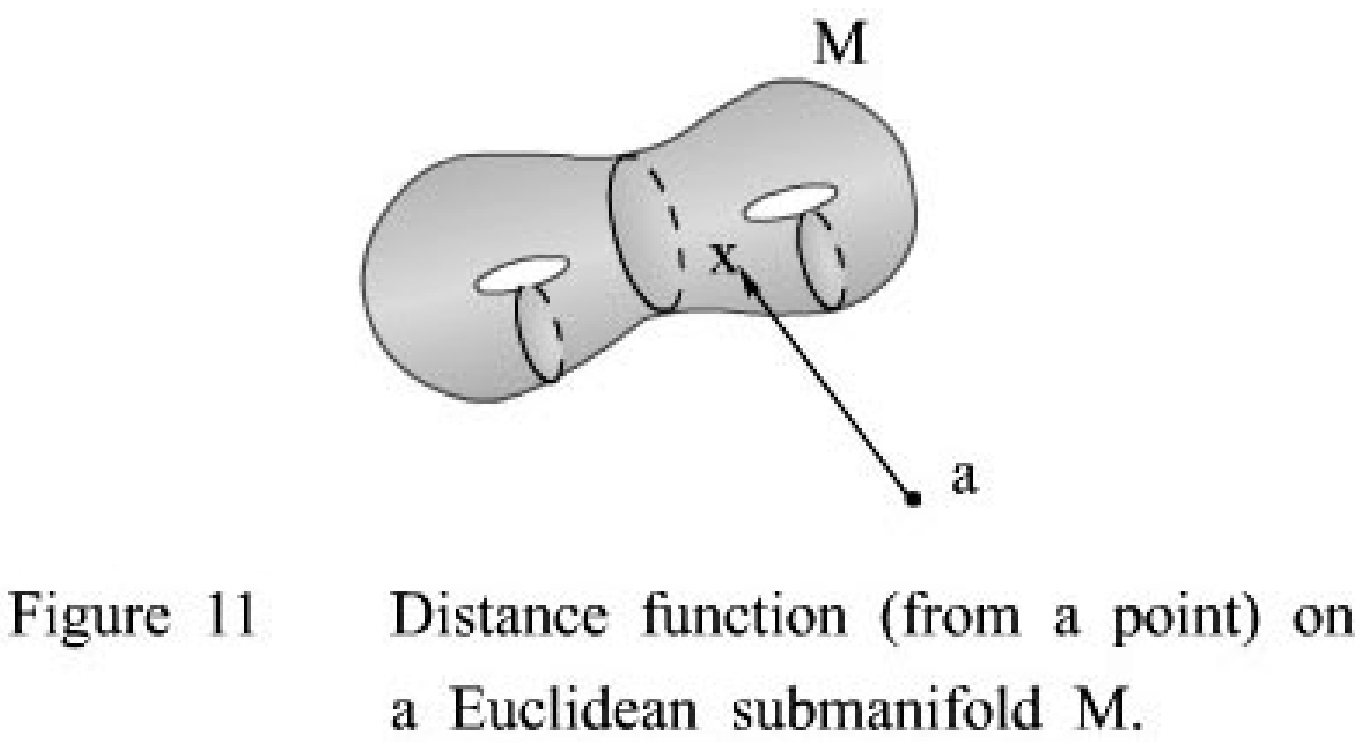}%
}%

\bigskip

\noindent Let $\Sigma_{a}$ be the set of all critical points of this function.
Two questions are:

\begin{quote}
(a) How to specify the critical set of $f_{a}$?

(b) For which choice of the point $a\in E$, $f_{a}$ is a Morse function on $M$?
\end{quote}

\bigskip

For a point $x\in M$ let $T_{x}M\subset E$ be the tangent plane to $M$ at $x$
(an affine plane in $E$ with dimension $n$). Its orthonormal complement

\begin{center}
$\gamma_{x}=\{v\in E\mid v\perp M_{x}\}$
\end{center}

\noindent is called \textsl{the normal plane to }$M$\textsl{ at }$x$. We state
the answers to questions (a) and (b) in

\textbf{Lemma 3.1.} \textsl{Let }$f_{a}:M\rightarrow\mathbb{R}$\textsl{ be as above.}

\textsl{(1) }$\Sigma_{a}=\{x\in M\mid a-x\in\gamma_{x}\}$\textsl{;}

\textsl{(2) For almost all }$a\in E$\textsl{, }$f_{a}$\textsl{ is a Morse function.}

\textbf{Proof.} The function $g_{a}:E\rightarrow\mathbb{R}$ by $x\rightarrow
\parallel x-a\parallel^{2}$ has gradient field grad$_{x}$ $g_{a}=2(x-a)$.
Since $f_{a}=g_{a}\mid M$, for a $x\in M$,

\begin{center}
grad$_{x}$ $f_{a}=$the orthonormal projection of $2(x-a)$ to $T_{x}M$.
\end{center}

\noindent So $x\in\Sigma_{a}$ (i.e. grad$_{x}$ $f_{a}=0$) is equivalent to
$2(x-a)\perp T_{x}M$. This shows (1).

Let $\Lambda\subset E$ be the focal set of the submanifold $M\subset E$. It
can be shown that $f_{a}$ is a Morse function if and only if $a\in
E\backslash\Lambda$. (2) follows from the fact that $\Lambda$ has measure $0$
in $E$ (cf. [M$_{2}$, p.32-38]). $\square$

\bigskip

\bigskip

\begin{center}
\textbf{3--2. Examples of submanifolds in Euclidean spaces}
\end{center}

Many manifolds important in geometry are already sitting in Euclidean spaces
in some ready-made fashion. We present such examples.

Let $\mathbb{F}$ be one of $\mathbb{R}$ (the field of reals)$,\mathbb{C}$ (the
field of complex) or $\mathbb{H}$ (the division algebra of quaternions). Let
$E$ be one of the following real vector spaces:

the space of $n\times n$ matrices over $\mathbb{F}$: $M(n;\mathbb{F})$;

the space of complex Hermitian matrices:

\begin{center}
$S(n;\mathbb{C})=\{x\in M(n;\mathbb{C})\mid x^{\tau}=x\}$;
\end{center}

the space of complex symmetric matrices

\begin{center}
$S^{+}(n;\mathbb{C})=\{x\in M(n;\mathbb{C})\mid x^{\tau}=\overline{x}\}$;
\end{center}

the space of real skew symmetric matrices:

\begin{center}
$S^{-}(2n;\mathbb{R})=\{x\in M(2n;\mathbb{R})\mid x^{\tau}=-x\}$.
\end{center}

\noindent Their dimensions as real vector spaces are respectively

$\qquad\dim_{\mathbb{R}}M(n;\mathbb{F})=\dim_{\mathbb{R}}\mathbb{F}\cdot
n^{2}$;

$\qquad\dim_{\mathbb{R}}S(n;\mathbb{C})=n(n+1)$;

$\qquad\dim_{\mathbb{R}}S^{+}(n;\mathbb{C})=n(n-1)$;

$\qquad\dim_{\mathbb{R}}S^{-}(2n;\mathbb{R})=n(2n-1)$.

\noindent Further, $E$ is an Euclidean space with the metric specified by

\begin{center}
$<x,y>=\operatorname{Re}[Tr(x^{\ast}y)]$, $x,y\in E$,
\end{center}

\noindent where $\ast$ means transpose followed by conjugation.

\bigskip

Consider in $E$ the following submanifolds

$\qquad O(n;\mathbb{F})=\{x\in M(n;\mathbb{F})\mid x^{\ast}x=I_{n}\}$

$\qquad G_{n,k}=\{x\in S^{+}(n;\mathbb{C})\mid x^{2}=I_{n}$, $l(x)=k\}$;

$\qquad LG_{n}=\{x\in S(n;\mathbb{C})\mid\overline{x}x=I_{n}\}$;

\qquad$\mathbb{C}S_{n}=\{x\in S^{-}(2n;\mathbb{R})\mid x^{2}=-I_{2n}\}$,

\noindent where $l(x)$ means ``the number of negative eigenvalues of $x$''and
where $I_{n}$ is the identity matrix. The geometric interests in these
manifolds may be illustrated in

\begin{center}
$O(n;\mathbb{F})=\{%
\begin{array}
[c]{c}%
O(n)\text{ if }\mathbb{F}=\mathbb{R}\text{: the real orthogonal group of rank
}n\text{;~}\\
U(n)\text{ if }\mathbb{F}=\mathbb{C}\text{: the unitary group of rank
}n\text{;\qquad\qquad}\\
Sp(n)\text{ if }\mathbb{F}=\mathbb{H}\text{: the symplectic group of rank
}n\text{;~\qquad}%
\end{array}
$
\end{center}

$G_{n,k}$: the Grassmannian of $k$-subspaces in $\mathbb{C}^{n}$;

$LG_{n}$: the Grassmannian of Largrangian subspaces in $\mathbb{C}^{n}$;

$\mathbb{C}S_{n}$: the Grassmannian of complex structures on $\mathbb{R}^{2n}$.

\begin{center}
\textbf{3--3.Morse functions via Euclidean geometry}
\end{center}

Let $0<\lambda_{1}<\cdots<\lambda_{n}$ be a sequence of $n$ reals, and let
$a\in E$ be the point with

\begin{center}
$a=\{%
\begin{array}
[c]{c}%
diag\{\lambda_{1},\cdots,\lambda_{n}\}\text{ if }M\neq\mathbb{C}%
S_{n}\text{;\qquad\qquad\qquad\qquad\qquad\qquad}\\
\lambda_{1}J\oplus\cdots\oplus\lambda_{n}J\text{, }J=\left(
\begin{array}
[c]{cc}%
0 & -1\\
1 & 0
\end{array}
\right)  \text{, if }M=\mathbb{C}S_{n}\text{.\quad\qquad\quad}%
\end{array}
$
\end{center}

\noindent With respect to the metric on $E$ specified in 3-2, the function

\begin{center}
$f_{a}:M\rightarrow\mathbb{R}$, $f_{a}(x)=\parallel x-a\parallel^{2}$
\end{center}

\noindent admits a simple-looking expression

\begin{center}
$f_{a}((x_{ij}))=<x,x>+<a,a>-2<a,x>$

\bigskip

$=const-2\{%
\begin{array}
[c]{c}%
\Sigma\lambda_{i}\operatorname{Re}(x_{ii})\text{ if }M=G_{n,k}\text{,
}O(n;\mathbb{F})\text{, }LG_{n}\text{; and}\\
\Sigma\lambda_{i}x_{2i-1,2i}\text{ if }M=\mathbb{C}S_{n}\text{.\qquad
\qquad\qquad\qquad}%
\end{array}
$
\end{center}

\noindent For a subsequence $I=[i_{1},\cdots,i_{r}]\subseteq\lbrack
1,\cdots,n]$, denote by $\sigma_{I}\in E$ the point

\begin{center}
$\sigma_{I}=\{%
\begin{array}
[c]{c}%
diag\{\varepsilon_{1},\cdots,\varepsilon_{n}\}\text{ if }M\neq\mathbb{C}%
S_{n}\text{;}\\
\varepsilon_{1}J\oplus\cdots\oplus\varepsilon_{n}J\text{ if }M=\mathbb{C}%
S_{n}\text{,}%
\end{array}
$
\end{center}

\noindent where $\varepsilon_{k}=-1$ if $k\in I$ and $\varepsilon_{k}=1$ otherwise.

\bigskip

\textbf{Theorem 3.} \textsl{In each of the above four cases, }$f_{a}%
:M\rightarrow R$\textsl{ is a Morse function on }$M$\textsl{. Further,}

\textsl{(1) the set of critical points of }$f_{a}$\textsl{ is}

\begin{center}
$\Sigma_{a}=\{%
\begin{array}
[c]{c}%
\{\sigma_{0},\text{ }\sigma_{I}\in M\mid I\subseteq\lbrack1,\cdots,n]\}\text{
if }M\neq G_{n,k}\text{;\qquad\quad\quad}\\
\{\sigma_{I}\in M\mid I\subseteq\lbrack1,\cdots,n]\text{ with}\mid
I\mid=k\}\text{ if }M=G_{n,k}\text{.}%
\end{array}
$
\end{center}

\textsl{(2) the index functions are given respectively by}

\begin{center}
$Ind(\sigma_{i_{1},\cdots,i_{r}})=\{%
\begin{array}
[c]{c}%
\dim_{\mathbb{R}}\mathbb{F}\cdot(i_{1}+\cdots+i_{r})-r\text{ if }%
M=O(n;\mathbb{F})\text{;\ }\\
2(i_{1}+\cdots+i_{r}-r)\text{ if }M=\mathbb{C}S_{n}\text{;\quad\qquad\qquad}\\
i_{1}+\cdots+i_{r}\text{ if }M=LG_{n}\text{;}\qquad\qquad\qquad\qquad
\end{array}
$

\bigskip

$Ind(\sigma_{i_{1},\cdots,i_{k}})=2\underset{1\leq s\leq k}{\Sigma}(i_{s}%
-s)$\textsl{ if }$M=G_{n,k}$\textsl{.}
\end{center}

\bigskip

\begin{center}
\textbf{3--4. Proof of Theorem 3}
\end{center}

We conclude Section 3 by a proof of Theorem 3.

\textbf{Lemma 3.2.} \textsl{For a }$x\in M$\textsl{ one has}

\begin{center}
$T_{x}M=\{%
\begin{array}
[c]{c}%
\{u\in E\mid xu=-ux\}\text{ for }M=G_{n,k}\text{; }\mathbb{C}S_{n}\\
\{u\in E\mid x^{\ast}u=-u^{\ast}x\}\text{ for }M=O(n;\mathbb{F})\\
\{u\in E\mid\overline{x}u=-\overline{u}x\}\text{ for }M=LG_{n}\text{.\qquad}%
\end{array}
$
\end{center}

\noindent\textsl{Consequently}

\begin{center}
$\gamma_{x}M=\{%
\begin{array}
[c]{c}%
\{u\in E\mid xu=ux\}\text{ for }M=G_{n,k}\text{; }\mathbb{C}S_{n}\\
\{u\in E\mid x^{\ast}u=u^{\ast}x\}\text{ for }M=O(n;\mathbb{F})\\
\{u\in E\mid\overline{x}u=\overline{u}x\}\text{ for }M=LG_{n}\text{.\qquad}%
\end{array}
$
\end{center}

\textbf{Proof. }We verify Lemma 3.2 for the case\textbf{ }$M=G_{n,k}$ as an
example. Consider the map $h:S^{+}(n;\mathbb{C})\rightarrow S^{+}%
(n;\mathbb{C})$ by $x\rightarrow x^{2}$. Then

\begin{quote}
(1) $h^{-1}(I_{n})=\underset{1\leq t\leq n-1}{\sqcup}G_{n,t}$;

(2) the tangent map of $h$ at a point $x\in S^{+}(n;\mathbb{C})$ is
\end{quote}

\begin{center}
$T_{x}h(u)=\underset{t\rightarrow0}{\lim}\frac{h(x+tu)-h(x)}{t}=ux+xu$.
\end{center}

\noindent It follows that, for a $x\in G_{n,k}$,

\begin{center}
$T_{x}G_{n,k}\subseteq KerT_{x}h=\{u\in S^{+}(n;\mathbb{C})\mid ux+xu=0\}$.
\end{center}

\noindent On the other hand $\dim_{\mathbb{C}}KerT_{x}h=k(n-k)$ ($=\dim
_{\mathbb{C}}T_{x}G_{n,k}$). So the dimension comparison yields

\begin{center}
$T_{x}G_{n,k}=\{u\in S^{+}(n;\mathbb{C})\mid xu=-ux\}$.
\end{center}

For any $x\in G_{n,k}$ the ambient space $E=S^{+}(n;\mathbb{C})$ admits the
orthogonal decomposition

\begin{center}
$S^{+}(n;\mathbb{C})=\{u\mid xu=-ux\}\oplus\{u\mid xu=ux\}$
\end{center}

\noindent in which the first summand has been identified with $T_{x}G_{n,k}$
in the above computation. It follows that $\gamma_{x}G_{n,k}=\{u\mid xu=ux\}$.

The other cases can be verified by the same method.$\square$

\bigskip

\textbf{Lemma 3.3.} \textsl{Statement (1) of Theorem 3 holds true.}

\textbf{Proof.} Consider the case $G_{n,k}\subset S^{+}(n;\mathbb{C})$.

$\ \ \ x\in\Sigma_{a}\Leftrightarrow x-a\in\gamma_{x}G_{n,k}$ (by (1) of Lemma 3.1)

$\qquad\qquad\Leftrightarrow(x-a)x=x(x-a)$ (by Lemma 3.2)

$\qquad\quad\quad\Leftrightarrow xa=ax$.

\noindent Since $a$ is diagonal with the distinguished diagonal entries
$\lambda_{1}<\cdots<\lambda_{n}$, $x$ is also diagonal. Since $x^{2}=I_{n}$
with $l(x)=k$, we must have $x=\sigma_{i_{1},\cdots,i_{k}}$ for some
$[i_{1},\cdots,i_{k}]\subseteq\lbrack1,\cdots,n]$.

Analogous computations verify the other cases.$\square$

\bigskip

To prove Theorem 3 we need examining the Hessian operator $H_{x_{0}}%
(f_{a}):T_{x_{0}}M\rightarrow T_{x_{0}}M$ \ at a critical point $x_{0}%
\in\Sigma_{a}$. The following formulae will be useful for this purpose.

\textbf{Lemma 3.4.} $H_{x_{0}}(f_{a})(u)=\{%
\begin{array}
[c]{c}%
(ua-au)x_{0}\text{ \ \ for }M=G_{n,k}\text{; }\mathbb{C}S_{n}\text{;}\\
(u^{\ast}a-au^{\ast})x_{0}\text{ for }M=O(n;\mathbb{F})\text{;\quad}\\
(\overline{u}a-a\overline{u})x_{0}\text{ for }M=LG_{n}\text{.~\qquad}%
\end{array}
$

\textbf{Proof.} As a function on the Euclidean space $E$, $f_{a}$ has gradient
field $2(x-a)$. However, the gradient field of the restricted function
$f_{a}\mid M$ \ is the orthogonal projection of $2(x-a)$ to $T_{x}M$.

In general, for any $x\in M$, a vector $u\in E$ has the ``canonical'' decomposition

\begin{center}
$u=\{%
\begin{array}
[c]{c}%
\frac{u-xux}{2}+\frac{u+xux}{2}\text{ if }M=G_{n,k}\text{; }\mathbb{C}%
S_{n}\text{;}\\
\frac{u-x^{\ast}ux}{2}+\frac{u+x^{\ast}ux}{2}\text{ if }M=O(n;\mathbb{F}%
)\text{;}\\
\frac{u-\overline{x}ux}{2}+\frac{u+\overline{x}ux}{2}\text{ if }%
M=LG_{n}\text{.\qquad}%
\end{array}
$
\end{center}

\noindent with the first component in the $T_{x}M$ and the second component in
$\gamma_{x}M$ by Lemma 3.2. Applying these to $u=2(x-a)$ yields respectively that

\begin{center}
$grad_{x}f_{a}=\{%
\begin{array}
[c]{c}%
(xax-a)\text{ for }M=G_{n,k}\text{; }\mathbb{C}S_{n}\text{;}\\
(x^{\ast}ax-a)\text{ for }M=O(n;\mathbb{F})\text{;\ }\\
(\overline{x}ax-a)\text{ for }M=LG_{n}\text{.\quad\quad}%
\end{array}
$
\end{center}

Finally, the Hessian operator can be computed in term of the gradient as

\begin{center}
$H_{x_{0}}(f_{a})(u)=\underset{t\rightarrow0}{\lim}\frac{grad_{x_{0}+tu}%
f_{a}-grad_{x_{0}}f_{a}}{t}$, $u\in T_{x}M$.
\end{center}

\noindent As an example we consider the case $M=G_{n,k}$. We have

$\underset{t\rightarrow0}{\qquad\qquad\lim}\frac{grad_{x_{0}+tu}%
f_{a}-grad_{x_{0}}f_{a}}{t}=$ $\underset{t\rightarrow0}{\lim}\frac
{[(x_{0}+tu)a(x_{0}+tu)-a]-[x_{0}ax_{0}-a]}{t}$

$=uax_{0}+x_{0}au=uax_{0}+ax_{0}u$ (because $a$ and $x_{0}$ are diagonal)

$=(ua-au)x_{0}$

\noindent(because vectors in $T_{x_{0}}G_{n,k}$ anti-commute with $x_{0}$ by
Lemma 3.2).$\square$

\bigskip

\textbf{Proof of Theorem 3.} In view of Lemma 3.3, Theorem 3 will be completed
once we have shown

\begin{quote}
(a) $f_{a}$ is non-degenerate at any $x_{0}\in\Sigma_{a}$; and

(b) the index functions on $\Sigma_{a}$ is given as that in (2) of Theorem 3.
\end{quote}

\noindent This can be done by applying Lemma 3.2 and Lemma 3.4. We verify
these for the cases $M=G_{n,k}$, $O(n)$ and $LG_{n}$ in detail, and leave the
other cases to the reader.

\textbf{Case 1.} $M=G_{n,k}\subset S^{+}(n;\mathbb{C})$.

\begin{quote}
(1) The most convenient vectors that span the real vector space $S^{+}%
(n;\mathbb{C})$ are
\end{quote}

\begin{center}
$\{b_{s,t}\mid1\leq s,t\leq n\}\sqcup\{c_{s,t}\mid1\leq s\neq t\leq n\}$,
\end{center}

\begin{quote}
\noindent where $b_{s,t}$ has the entry $1$ at the places $(s,t)$, $(t,s)$ and
$0$ otherwise, and where $c_{s,t}$ has the pure imaginary $i$ at $(s,t)$, $-i$
at the $(t,s)$ and $0$ otherwise.

(2) For a $x_{0}=\sigma_{I}\in\Sigma_{a}$, those $b_{s,t}$, $c_{s,t}$ that
``\textsl{anti-commute}'' with $x_{0}$ belong to $T_{x_{0}}G_{n,k}$ by Lemma
3.2, and form a basis for $T_{x_{0}}G_{n,k}$
\end{quote}

\begin{center}
$T_{x_{0}}G_{n,k}=\{b_{s,t},c_{s,t}\mid(s,t)\in I\times J\},$
\end{center}

\begin{quote}
\noindent where $J$ is the complement of $I$ in $[1,\cdots,n]$.

(3) Applying the Hessian (Lemma 3.4) to the $b_{s,t},c_{s,t}\in T_{x_{0}%
}G_{n,k}$ yields
\end{quote}

\begin{center}
$H_{x_{0}}(f_{a})(b_{s,t})=(\lambda_{t}-\lambda_{s})b_{s,t}$;

$H_{x_{0}}(f_{a})(c_{s,t})=(\lambda_{t}-\lambda_{s})c_{s,t}$.
\end{center}

\begin{quote}
\noindent That is, the $b_{s,t},c_{s,t}\in T_{x_{0}}G_{n,k}$ are precisely the
eigenvectors for the operator $H_{x_{0}}(f_{a})$. These indicate that
$H_{x_{0}}(f_{a})$ is nondegenerate (since $\lambda_{t}\neq\lambda_{s}$ for
all $s\neq t$), hence $f_{a}$ is a Morse function.

(4) It follows from the formulas in (3) that the negative space for $H_{x_{0}%
}(f_{a})$ is spanned by $\{b_{s,t}$, $c_{s,t}\mid(s,t)\in I\times J,t<s\}$. Consequently
\end{quote}

\begin{center}
$Ind(\sigma_{I})=2\#\{(s,t)\in I\times J\mid t<s\}=2\underset{1\leq s\leq
k}{\Sigma}(i_{s}-s)$.
\end{center}

\bigskip

\textbf{Case 2.} $M=O(n)\subset M(n;\mathbb{R})$.

\begin{quote}
(1) A natural set of vectors that spans the space $M(n;\mathbb{R})$ is
\end{quote}

\begin{center}
$\{b_{s,t}\mid1\leq s\leq t\leq n\}\sqcup\{\beta_{s,t}\mid1\leq s<t\leq n\}$,
\end{center}

\begin{quote}
\noindent where $b_{s,t}$ is as case 1, and where $\beta_{s,t}$ is the skew
symmetric matrix with entry $1$ at the $(s,t)$ place, $-1$ at the $(t,s)$
place and $0$ otherwise;

(2) For a $x_{0}=\sigma_{I}\in\Sigma_{a}$ those $b_{s,t}$, $\beta_{s,t}$
that``anti-commute'' with $x_{0}$ yields precisely a basis for
\end{quote}

\begin{center}
\noindent$T_{x_{0}}O(n)=\{\beta_{s,t}\mid(s,t)\in I\times I,J\times
J,s<t\}\sqcup\{b_{s,t}\mid(s,t)\in I\times J\}$
\end{center}

\begin{quote}
\noindent by Lemma 3.2, where $J$ is the complement of $I$ in $[1,\cdots,n]$.

(3) Applying the Hessian operator (Lemma 3.4) to $b_{s,t}$, $\beta_{s,t}\in
T_{x_{0}}O(n)$ tells
\end{quote}

\begin{center}
$H_{x_{0}}(f_{a})(\beta_{s,t})=\{%
\begin{array}
[c]{c}%
-(\lambda_{t}+\lambda_{s})\beta_{s,t}\text{ if }(s,t)\in I\times
I,s<t\text{;}\\
(\lambda_{t}+\lambda_{s})\beta_{s,t}\text{ if }(s,t)\in J\times J,\text{
}s<t\text{.}%
\end{array}
$

$H_{x_{0}}(f_{a})(b_{s,t})=(\lambda_{t}-\lambda_{s})b_{s,t}$ if $(s,t)\in
I\times J$.
\end{center}

\begin{quote}
\noindent This implies that the $b_{s,t}$, $\beta_{s,t}\in T_{x_{0}}G_{n,k}$
are precisely the eigenvectors for the operator $H_{x_{0}}(f_{a})$, and the
$f_{a\text{ }}$is a Morse function.

(4) It follows from the computation in (3) that

$Ind(\sigma_{I})=\#\{(s,t)\in I\times I\mid s<t\}+\#\{(s,t)\in I\times J\mid t<s\}$
\end{quote}

$=1+2+\cdots+(r-1)+[(i_{1}-1)+(i_{2}-2)+\cdots+(i_{r}-r)]$

$=\Sigma i_{s}-r$.

\bigskip

\textbf{Case 3.} $M=LG_{n}\subset S(n;\mathbb{C})$.

\begin{quote}
(1) Over reals, the most natural vectors that span the space $S(n;C)$ are
\end{quote}

\begin{center}
$\{b_{s,t}\mid1\leq s,t\leq n\}\cup\{ib_{s,t}\mid1\leq s,t\leq n\}$,
\end{center}

\begin{quote}
\noindent where $b_{s,t}$ is as that in Case 1 and where $i$ is the pure imaginary;

(2) For a $x_{0}=\sigma_{I}\in\Sigma_{a}$ those ``anti-commute'' with $x_{0}$
yields precisely a basis for $T_{x_{0}}LG_{n}$
\end{quote}

\begin{center}
$T_{x_{0}}LG_{n}=\{b_{s,t}\mid(s,t)\in I\times J\amalg J\times I\}\sqcup$

$\{ib_{s,t}\mid(s,t)\in I\times I\sqcup J\times J\}$
\end{center}

\begin{quote}
\noindent where $J$ is the complement of $I$ in $[1,\cdots,n]$.

(3) Applying the Hessian to $b_{s,t}$, $ib_{s,t}\in T_{x_{0}}LG_{n}$ (cf.
Lemma 3.4) tells
\end{quote}

\begin{center}
$H_{x_{0}}(f_{a})(ib_{s,t})=\{%
\begin{array}
[c]{c}%
-(\lambda_{t}+\lambda_{s})ib_{s,t}\text{ if }(s,t)\in I\times I\\
(\lambda_{t}+\lambda_{s})ib_{s,t}\text{ if }(s,t)\in J\times J
\end{array}
$;

$H_{x_{0}}(f_{a})(b_{s,t})=\{%
\begin{array}
[c]{c}%
(\lambda_{t}-\lambda_{s})b_{s,t}\text{ if }(s,t)\in I\times J\\
(\lambda_{s}-\lambda_{t})b_{s,t}\text{ if }(s,t)\in J\times I
\end{array}
$.
\end{center}

\begin{quote}
\noindent It follows that the $b_{s,t}$, $ib_{s,t}\in T_{x_{0}}G_{n,k}$ are
precisely the eigenvectors for the operator $H_{x_{0}}(f_{a})$, and $f_{a}$ is
a Morse function.

(4) It follows from (2) and (3) that

$Ind(\sigma_{I})=\#\{(s,t)\in I\times I\mid t\leq s\}+\#\{(s,t)\in I\times
J\mid t\leq s\}$
\end{quote}

$\qquad\qquad=i_{1}+\cdots+i_{r}$.$\square$

\bigskip

\textbf{Remark.} Let $E$ be one of the following matrix spaces:

\qquad the space of $n\times k$ matrices over $\mathbb{F}$: $M(n\times
k;\mathbb{F})$;

\qquad the space of symmetric matrices $S^{+}(n;\mathbb{F})=\{x\in
M(n;\mathbb{F})\mid x^{\tau}=\overline{x}\}$.

\noindent Consider in $E$ the following submanifolds:

\begin{center}
$V_{n,k}(\mathbb{F})=\{x\in M(n\times k;\mathbb{F})\mid\overline{x}^{\tau
}x=I_{k}\}$;

$G_{n,k}(\mathbb{F})=\{x\in S^{+}(n;\mathbb{F})\mid x^{2}=I_{n}$, $l(x)=k\}$.
\end{center}

\noindent These are known respectively as the \textsl{Stiefel manifold} of
orthonormal $k$-frames on $\mathbb{F}^{n}$ (the $n$-dimensional $\mathbb{F}%
$-vector space) and the \textsl{Grassmannian} of $k$-dimen-

\noindent sional $\mathbb{F}$-subspaces in $\mathbb{F}^{n}$. Results analogous
to Theorem 3 hold for these two family of manifolds as well [D$_{1}$],
[D$_{2}$].

\textbf{Remark. }In [VD, Theorem 1.2], the authors proved that the function
$f_{a}$ on $M=G_{n,k}(\mathbb{F})$, $LG_{n}$, $\mathbb{C}S_{n}$ is a perfect
Morse function (without specifying the set $\Sigma_{a}$ as well as the index
function Ind: $\Sigma_{a}\rightarrow\mathbb{Z}$).

\section{Morse functions of Bott-Samelson type}

We recall the original construction of Bott-Samelson cycles in 4--1 and
explain its generalization due to Hsiang-Palais-Terng [HTP] in 4--2.

In fact, the Morse functions concerned in Theorem 3 are all Bott-Samelson type
(cf. Theorem 6). The induced cohomology homomorphism of Bott-Samelson cycles
enables one to resolve the multiplication in cohomology into the
multiplication of symmetric functions of various types (Theorem 7).

\begin{center}
\textbf{4--1. Morse functions on flag manifolds }(cf. [BS$_{1}$,BS$_{2}$]).
\end{center}

Let $G$ be a compact connected semi-simple Lie group with the unit $e\in G$
and a fixed maximal torus $T\subset G$. The tangent space $L(G)=T_{e}G$ (resp.
$L(T)=T_{e}T$) is canonically furnished with the structure of an algebra,
known as the \textsl{Lie algebra} (resp. the \textsl{Cartan subalgebra}) of
$G$. The exponential map induces the commutative diagram

\begin{center}
$%
\begin{array}
[c]{ccc}%
L(T) & \rightarrow & L(G)\\
\qquad\qquad\exp\downarrow\qquad\qquad &  & \qquad\downarrow\exp\\
T & \rightarrow & G
\end{array}
$
\end{center}

\noindent where the horizontal maps are the obvious inclusions. Equip $L(G)$
(hence also $L(T)$) an inner product invariant under the adjoint action of $G$
on $L(G)$.

For a $v\in L(T)$ let $C(v)$ be the centralizer of $\exp(v)\in G$. The set of
singular points in $L(T)$ is the subspace of the Cartan subalgebra $L(T)$:

\begin{center}
$\Gamma=\{v\in L(T)\mid\dim C(v)>\dim T\}$.
\end{center}

\textbf{Lemma 4.1.} \textsl{Let }$m=\frac{1}{2}(\dim G-\dim T)$\textsl{. There
are precisely }$m$\textsl{ hyperplanes }$L_{1},\cdots,L_{m}\subset
L(T)$\textsl{ through the origin }$0\in L(T)$\textsl{ so that }$\Gamma
=\underset{1\leq i\leq m}{\cup}L_{i}$\textsl{.}$\square$

$\qquad$

The planes $L_{1},\cdots,L_{m}$ are known as the singular planes of $G$. It
divide $L(T)$ into finite many convex hulls, known as the \textsl{Weyl
chambers} of $G$. Reflections in these planes generate the \textsl{Weyl group}
$W$ of $G$.

\bigskip

Fix a regular point $a\in L(T)$. The adjoint representation of $G$ gives rise
to a map $G\rightarrow L(G)$ by $g\rightarrow Ad_{g}(a)$, which induces an
embedding of the \textsl{flag manifold} $G/T=\{gT\mid g\in G\}$ of left cosets
of $T$ in $G$ into $L(G)$. In this way $G/T$ becomes a submanifold in the
Euclidean space $L(G)$.

Consider the function $f_{a}:G/T\rightarrow\mathbb{R}$ by $f_{a}(x)=\parallel
x-a\parallel^{2}$. The following beautiful result of Bott and Samelson
[BS$_{1}$,BS$_{2}$] tells how to read the critical points information of
$f_{a}$ from the linear geometry of the vector space $L(T)$.

\textbf{Theorem 4.} $f_{a}$\textsl{ is a Morse function on }$G/T$\textsl{ with
critical set}

\begin{center}
$\Sigma_{a}=\{w(a)\in L(T)\mid w\in W\}$
\end{center}

\noindent\textsl{(the orbit of the }$W$\textsl{-action on }$L(T)$\textsl{
through the point }$a\in L(T)$\textsl{).}

\textsl{The index function Ind}$:\Sigma_{a}\rightarrow\mathbb{Z}$\textsl{ is
given by}

\begin{center}
\textsl{Ind}$(w(a))=2\#\{L_{i}\mid L_{i}\cap\lbrack a,w(a)]\neq\emptyset]$\textsl{,}
\end{center}

\noindent\textsl{where }$[a,w(a)]$\textsl{ is the segment in }$L(T)$\textsl{
from }$a$\textsl{ to }$w(a)$\textsl{.}

\bigskip

Moreover, Bott and Samelson constructed a set of geometric cycles in $G/T$
that realizes an additive basis of $H_{\ast}(G/T;\mathbb{Z})$ as follows.

For a singular plane $L_{i}\subset L(T)$ let $K_{i}\subset G$ be the
centralizer of $\exp(L_{i})$. The Lie subgroup $K_{i}$ is \textsl{very simple}
in the sense that $T\subset K_{i}$ is also a maximal torus with the quotient
$K_{i}/T$ diffeomorphic to the $2$-sphere $S^{2}$.

For a $w\in W$ assume that the singular planes that meet the directed segments
$[a,w(a)]$ are in the order $L_{1},\cdots,L_{r}$. Put $\Gamma_{w}=K_{1}%
\times_{T}\cdots\times_{T}K_{r}$, where the action of $T\times\cdots\times T$
($r$-copies) acts on $K_{1}\times\cdots\times K_{r}$ from the left by

\begin{center}
$(k_{1},\cdots,k_{r})(t_{1},\cdots,t_{r})=(k_{1}t_{1},t_{1}^{-1}k_{2}%
t_{2},\cdots,t_{r-1}^{-1}k_{r}t_{r})$.
\end{center}

\noindent The map $K_{1}\times\cdots\times K_{r}\rightarrow G/T$ by

\begin{center}
$(k_{1},\cdots,k_{r})\rightarrow Ad_{k_{1}\cdots k_{r}}(w(a))$
\end{center}

\noindent clearly factors through the quotient manifold $\Gamma_{w}$, hence
induces a map

\begin{center}
$g_{w}:\Gamma_{w}\rightarrow G/T$.
\end{center}

\textbf{Theorem 5.} \textsl{The homology }$H_{\ast}(G/T;\mathbb{Z})$\textsl{
is torsion free with the additive basis }$\{g_{w\ast}[\Gamma_{w}]\in H_{\ast
}(G/T;\mathbb{Z})\mid w\in W\}$\textsl{.}

\textbf{Proof.} Let $e\in K_{i}$($\subset G$) be the group unit and put
$\overline{e}=[e,\cdots,e]\in\Gamma_{w}$. It were actually shown by Bott and
Samelson that

\begin{quote}
(1) $g_{w}^{-1}(w(a))$ consists of the single point $\overline{e}$;

(2) the composed function $f_{a}\circ g_{w}:\Gamma_{w}\rightarrow\mathbb{R}$
attains its maximum only at $\overline{e}$;

(3) the tangent map of $g_{w}$ at $\overline{e}$ maps the tangent space of
$\Gamma_{w}$ at $\overline{e}$ isomorphically onto the negative part of
$H_{w(a)}(f_{a})$.
\end{quote}

\noindent The proof is completed by Lemma 4.2 in \textbf{4.2}.$\square$

\bigskip

\textbf{Remark.} It was shown by Chevalley in 1958 [Ch] that the flag manifold
$G/T$ admits a cell decomposition $G/T=\underset{w\in W}{\cup}X_{w}$ indexed
by elements in $W$, with each cell $X_{w}$ an algebraic variety, known as a
Schubert variety on $G/T$. Hansen [Han] proved in 1971 that $g_{w}(\Gamma
_{w})=X_{w}$, $w\in W$. So the map $g_{w}$ is currently known as the
``\textsl{Bott-Samelson resolution} of $X_{w}$''.

For the description of Bott-Samelson cycles and their applications in
Algebro-geometric setting, see M. Brion [Br] in this volume.

\begin{center}
\textbf{4--2. Morse function of Bott-Samelson type}
\end{center}

In differential geometry, the study of \textsl{isoparametric submanifolds} was
begun by E. Cartan in 1933. In order to generalize Bott-Samelson's above cited
results to these manifolds Hsiang, Palais and Terng introduced the following
notation in their work [HPT]\footnote{In fact, the embedding $G/T\subset L(G)$
described in 4-1 defines $G/T$ as an isoparametric submanifold in $L(G)$
[HPT].}.

\ \textbf{Definition 4.1.} \textsl{A Morse function }$f:M\rightarrow
R$\textsl{ on a smooth closed manifold is said to be of Bott-Samelson type
over }$\mathbb{Z}_{2}$\textsl{(resp. }$\mathbb{Z}$\textsl{) if for each }%
$p\in\Sigma_{f}$\textsl{ there is a map (called a Bott-Samelson cycle of }%
$f$\textsl{ at }$p$\textsl{)}

\begin{center}
$g_{p}:N_{p}\rightarrow M$
\end{center}

\noindent\textsl{where }$N_{p}$\textsl{ is a closed oriented (resp.
unoriented) manifold of dimension }$Ind(p)$\textsl{ and where}

\begin{quote}
\textsl{(1) }$g_{p}^{-1}(p)=\{\overline{p}\}$\textsl{ (a single point);}

\textsl{(2) }$f\circ g_{p}$\textsl{ attains absolute maximum only at
}$\overline{p}$\textsl{;}

\textsl{(3) the tangent map }$T_{\overline{p}}g_{p}:T_{\overline{p}}%
N_{p}\rightarrow T_{p}M$\textsl{ is an isomorphism onto the negative space of
}$H_{p}(f)$\textsl{.}

\ \ \ \
{\includegraphics[
natheight=3.257800in,
natwidth=5.904900in,
height=2.0029in,
width=3.8363in
]%
{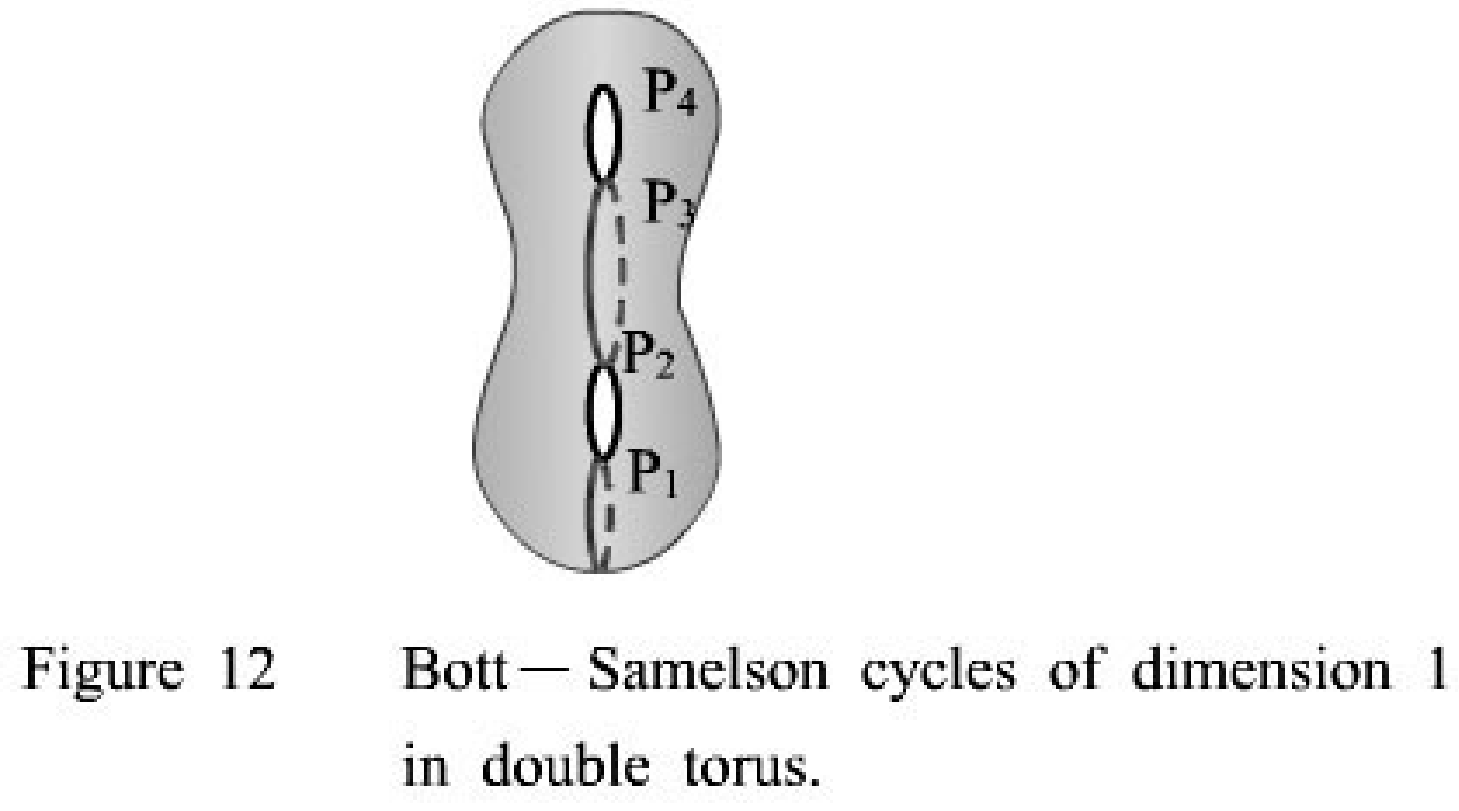}%
}%
\end{quote}

Information that one can get from a Morse function of Bott-Samelson type can
be seen from the next result [HPT].

\textbf{Lemma 4.2.} \textsl{If }$f:M\rightarrow R$\textsl{ is a Morse function
of Bott-Samelson type with Bott-Samelson cycles }$\{g_{p}:N_{p}\rightarrow
M\mid p\in\Sigma_{f}\}$\textsl{, then }$H_{\ast}(M;\mathbb{Z})$\textsl{ (resp.
}$H_{\ast}(M;\mathbb{Z}_{2})$\textsl{) has the additive basis}

\begin{center}
$\{g_{p\ast}[N_{p}]\in H_{\ast}(M;\mathbb{Z})\mid p\in\Sigma_{f}\}$

\textsl{(resp. }$\{g_{p\ast}[N_{p}]_{2}\in H_{\ast}(M;\mathbb{Z}_{2})\mid
p\in\Sigma_{f}\}$\textsl{),}
\end{center}

\noindent\textsl{where }$g_{p\ast}:H_{\ast}(N_{p};\mathbb{Z})\rightarrow
H_{\ast}(M;\mathbb{Z})$\textsl{ is the induced homomorphism and where }%
$[N_{p}]\in H_{\ast}(N_{p};\mathbb{Z})$\textsl{ (resp. }$[N_{p}]_{2}\in
H_{\ast}(N_{p};\mathbb{Z}_{2})$\textsl{) is the orientation class (resp.
}$\mathbb{Z}_{2}$\textsl{-orientation class).}

\textbf{Proof.} Without loss of generalities we may assume (as in the proof of
Theorem 2) that $\Sigma_{f}=\{p_{1},\cdots,p_{m}\}$ and that $f(p_{k}%
)<f(p_{k+1})$, $1\leq k\in m-1$. Consider the filtration on $M$: $M_{1}\subset
M_{2}\subset\cdots\subset M_{m}=M$ defined by $f$ and $\Sigma_{f}$ such that
$M_{k+1}\backslash M_{k}$ contains $p_{k}$ for every $1\leq k\leq m-1$.

It suffices to show, that if we put $p=p_{k+1}$, $m=Ind(p)$, then

\begin{enumerate}
\item[(D)] $\qquad\qquad\qquad H_{r}(M_{k+1};\mathbb{Z})=\{%
\begin{array}
[c]{c}%
H_{r}(M_{k};\mathbb{Z})\text{ if }r\neq m\text{;\qquad}\\
H_{r}(M_{k};\mathbb{Z})\oplus\mathbb{Z}\text{ if }r=m\text{,}%
\end{array}
$
\end{enumerate}

\noindent where the summand $\mathbb{Z}$ is generated by $g_{p\ast}[N_{p}]$.

The Bott-Samelson cycle $g_{p}:N_{p}\rightarrow M$ (cf. Definition 4.1) is
clearly a map into $M_{k+1}$. Let $r:M_{k+1}\rightarrow M_{k}\cup D^{m}$ be
the strong deformation retraction from the proof of Theorem 2, and consider
the composed map

\begin{center}
$g:N_{p}\overset{g_{p}}{\rightarrow}M_{k+1}\overset{r}{\rightarrow}M_{k}\cup
D^{m}$.
\end{center}

\noindent The geometric constraints (1)-(3) on the Bott-Samelson cycle $g_{p}$
imply \ that there exists an Euclidean neighborhood $U\subset D^{m}$ centered
at $p=0\in D^{m}$ so that if $V=:g^{-1}(U)$, then $g$ restricts to a
diffeomorphism $g\mid V:V\rightarrow U$. The proof of (D) (hence of Lemma 4.2)
is clearly done by the exact ladder induced by the ``\textsl{relative
homeomorphism}''$g:(N_{p},N_{p}\backslash V)\rightarrow(M_{k}\cup D^{m}%
,M_{k}\cup D^{m}\backslash U)$

\begin{center}
\ \ \ \ \ \ \ \ \ \ \ \ $\ \ \ \ \ \ \ \ \ \ \ \ \ \ \ \underset{\parallel
}{\mathbb{Z}}$\ \ \ \ \ $\ \ \ \ \ \ \ \ \ \ \ \underset{\parallel}%
{\mathbb{Z}}$

$\qquad\ \ \ \ \ \ \ \ \ \ \ \ \ 0\rightarrow H_{m}(N_{p})\overset{\cong
}{\rightarrow}H_{m}(N_{p},N_{p}\backslash V)\rightarrow H_{d-1}(N_{p}%
\backslash V)\rightarrow\cdots$

\ \ \ \ \ \ \ \ \ \ \ \ \ \ \ \ \ \ \ \ \ \ $\ \ g_{\ast}\downarrow$
\ \ \ \ \ \ \ \ \ \ \ $\ g_{\ast}\downarrow\cong$

\noindent$0\rightarrow H_{d}(M_{k})\rightarrow H_{d}((M_{k}\cup D^{m}%
)\rightarrow H_{d}((M_{k}\cup D^{m},M_{k})\rightarrow H_{d-1}(M_{k}%
)\rightarrow\cdots$
\end{center}

\noindent.$\square$

\begin{center}
\textbf{4--3. Bott-Samelson cycles and resolution of Schubert varieties}
\end{center}

Le $M$ be one of the following manifolds

$\qquad O(n;\mathbb{F})$: orthogonal (or unitary, or symplectic) group of rank
$n$;

$\qquad\mathbb{C}S_{n}$: the Grassmannian of complex structures on
$\mathbb{R}^{2n}$;

$\qquad G_{n,k}$: the Grassmannian of $k$-linear subspaces on $\mathbb{C}^{n}$

\noindent and

$\qquad LG_{n}$: the Grassmannian of Lagrangian subspaces on $\mathbb{C}^{n}$.

\noindent Let $f_{a}:M\rightarrow\mathbb{R}$ be the Morse function considered
in Theorem 3 of \S3.

\textbf{Theorem 6}. \textsl{In each case }$f_{a}$\textsl{ is a Morse function
of Bott-Samelson type which is}

\textsl{(1) over }$\mathbb{Z}$\textsl{ for }$M=U(n),Sp(n),$\textsl{
}$\mathbb{C}S_{n},G_{n,k}$\textsl{;}

\textsl{(2) over }$\mathbb{Z}_{2}$\textsl{ for }$M=O(n)$\textsl{ and }$LG_{n}$.

\bigskip

Instead of giving a proof of this result I would like to show the geometric
construction of the Bott-Samelson cycles required to justify the theorem, and
to point out the consequences which follow up (cf. Theorem 7).

Let $\mathbb{R}P^{n-1}$ be the real projective space of lines through the
origin $0$ in $\mathbb{R}^{n}$; $\mathbb{C}P^{n-1}$ the complex projective
space of complex lines through the origin $0$ in $\mathbb{C}^{n}$, and let
$G_{2}(\mathbb{R}^{2n})$ be the Grassmannian of oriented $2$-planes through
the origin in $\mathbb{R}^{2n}$.

\textbf{Construction 1}. Resolution $h:\widetilde{M}\rightarrow M$ of $M$.

\begin{quote}
(1) If $M=SO(n)$ (the special orthogonal group of order $n$) we let
$\widetilde{M}=\mathbb{R}P^{n-1}\times\cdots\times\mathbb{R}P^{n-1}$
($n^{\prime}$-copies, where $n^{\prime}=2[\frac{n}{2}]$) and define the map
$h:\widetilde{M}\rightarrow M$ to be
\end{quote}

\begin{center}
$h(l_{1},\cdots,l_{n^{\prime}})=\Pi_{1\leq i\leq n^{\prime}}R(l_{i})$,
\end{center}

\begin{quote}
\noindent where $l_{i}\in\mathbb{R}P^{n-1}$ and where $R(l_{i})$ is the
reflection on $\mathbb{R}^{n}$ in the hyperplane $l_{i}^{\perp}$ orthogonal to
$l_{i}$.

(2) If $M=G_{n,k}$ we let
\end{quote}

\begin{center}
$\widetilde{M}=\{(l_{1},\cdots,l_{k})\in\mathbb{C}P^{n-1}\times\cdots
\times\mathbb{C}P^{n-1}\mid l_{i}\perp l_{j}\}$ ($k$-copies)
\end{center}

\begin{quote}
\noindent and define the map $h:\widetilde{M}\rightarrow M$ to be
$h(l_{1},\cdots,l_{k})=<l_{1},\cdots,l_{k}>$, where $l_{i}\in\mathbb{C}%
P^{n-1}$ and where $<l_{1},\cdots,l_{k}>$ means the k-plane spanned by the
$l_{1},\cdots,l_{k}$.

(3) If $M=\mathbb{C}S_{n}$ we let
\end{quote}

\begin{center}
$\widetilde{M}=\{(L_{1},\cdots,L_{n})\in G_{2}(\mathbb{R}^{2n})\times
\cdots\times G_{2}(\mathbb{R}^{2n})\mid L_{i}\perp L_{j}\}$ ($n$-copies)
\end{center}

\begin{quote}
\noindent and define the map $h:\widetilde{M}\rightarrow M$ to be
$h(L_{1},\cdots,L_{k})=\Pi_{1\leq i\leq n}\tau(L_{i})$, where $L_{i}\in
G_{2}(\mathbb{R}^{2n})$ and where $\tau(L_{i}):\mathbb{R}^{2n}\rightarrow
\mathbb{R}^{2n}$ is the isometry which fixes points in the orthogonal
complements $L_{i}^{\perp}$ of $L_{i}$ and is the $\frac{\pi}{2}$ rotation on
$L_{i}$ in accordance with the orientation.
\end{quote}

\bigskip

\textbf{Construction 2}. Bott-Samelson cycles for the Morse function
$f_{a}:M\rightarrow\mathbb{R}$ (cf. [section3, Theorem 3]).

\begin{quote}
(1) If $M=SO(n)$ then $\Sigma_{a}=\{\sigma_{0},$ $\sigma_{I}\in M\mid
I\subseteq\lbrack1,\cdots,n],\mid I\mid\leq n^{\prime}\}$. For each
$I=(i_{1},\cdots,i_{r})\subseteq\lbrack1,\cdots,n]$ we put
\end{quote}

\begin{center}
$\mathbb{R}P[I]=$ $\mathbb{R}P^{0}\times\cdots\times\mathbb{R}P^{0}%
\times\mathbb{R}P^{i_{1}}\times\cdots\times\mathbb{R}P^{i_{r}}$ ($n^{\prime}$-copies).
\end{center}

\begin{quote}
\noindent Since $\mathbb{R}P[I]\subset\widetilde{M}$ we may set $h_{I}%
=h\mid\mathbb{R}P[I]$.
\end{quote}

\textsl{The map }$h_{I}:\mathbb{R}P[I]\rightarrow SO(n)$\textsl{ is a
Bott-Samelson cycle for }$f_{a}$\textsl{ at }$\sigma_{I}$.

\begin{quote}
(2) If $M=G_{n,k}$ then $\Sigma_{a}=\{\sigma_{I}\in M\mid I=(i_{1}%
,\cdots,i_{k})\subseteq\lbrack1,\cdots,n]\}$. For each $I=(i_{1},\cdots
,i_{k})\subseteq\lbrack1,\cdots,n]$ we have
\end{quote}

\begin{center}
$\mathbb{C}P^{i_{1}}\times\cdots\times\mathbb{C}P^{i_{k}},\widetilde{M}%
\subset\mathbb{C}P^{n-1}\times\cdots\times\mathbb{C}P^{n-1}$($k$-copies).
\end{center}

\begin{quote}
\noindent So we may define the intersection $\mathbb{C}P[I]=$ $\mathbb{C}%
P^{i_{1}}\times\cdots\times\mathbb{C}P^{i_{k}}\cap\widetilde{M}$ in
$\mathbb{C}P^{n-1}\times\cdots\times\mathbb{C}P^{n-1}$ and set $h_{I}%
=h\mid\mathbb{C}P[I]$.
\end{quote}

\textsl{The map }$h_{I}:\mathbb{C}P[I]\rightarrow G_{n,k}$\textsl{ is a
Bott-Samelson cycle for }$f_{a}$\textsl{ at }$\sigma_{I}$.

\bigskip

\begin{center}
\textbf{4--4. Multiplication in cohomology: }\textsl{Geometry versus combinatorics}
\end{center}

Up to now we have plenty examples of Morse functions of Bott-Samelson type.
Let $f:M\rightarrow\mathbb{R}$ be such a function with critical set
$\Sigma_{f}=\{p_{1},\cdots,p_{m}\}$. From the proof of Lemma 4.2 we see that
each descending cell $S(p_{i})\subset M$ forms a closed cycle on $M$ and all
of them form an additive basis for the homology

\begin{center}
$\{[S(p_{i})]\in H_{r_{i}}(M;\mathbb{Z}$ or $\mathbb{Z}_{2}$)$\mid1\leq i\leq
m$, $r_{i}=Ind(p_{i})\}$,
\end{center}

\noindent where the coefficients in homology depend on whether the
Bott-Samelson cycles are orientable or not.

Many pervious work on Morse functions stopped at this stage, for people were
content to have found Morse functions on manifolds whose critical points
determine an additive basis for homology (such functions are normally called
\textsl{perfect Morse functions}).

\bigskip

However, the difficult task that one has experienced in topology is not to
find an additive basis for homology, but is to understand the multiplicative
rule among basis elements in cohomology. More precisely, we let

\begin{center}
$\{[\Omega(p_{i})]\in H^{r_{i}}(M;\mathbb{Z}$ or $\mathbb{Z}_{2})\mid1\leq
i\leq m$, $r_{i}=Ind(p_{i})\}$
\end{center}

\noindent be the basis for the cohomology Kronecker dual to the $[S(p_{i})]$ as

\begin{center}
$<[\Omega(p_{i})],[S(p_{j})]>=\delta_{ij}$.
\end{center}

\noindent Then we must have the expression

\begin{center}
$[\Omega(p_{i})]\cdot\lbrack\Omega(p_{j})]=\Sigma a_{ij}^{k}[\Omega(p_{k})]$
\end{center}

\noindent in the ring $H^{\ast}(M;\mathbb{Z}$ or $\mathbb{Z}_{2})$, where
$a_{ij}^{k}\in\mathbb{Z}$ or $\mathbb{Z}_{2}$ depending on whether the
Bott-Samelson cycles orientable or not, and where $\cdot$ means
\textsl{intersection product} in Algebraic Geometry and \textsl{cup product}
in Topology.

\textbf{Problem 4}. \textsl{Find the numbers }$a_{ij}^{k}$\textsl{ for each
triple }$1\leq i,j,k\leq m$\textsl{.}

\bigskip

To emphasis Problem 4 we quote from N. Steenrod [St, p.98]:

\begin{quote}
``the cup product requires a diagonal approximation $d_{\#}:M\rightarrow
M\times M$. Many difficulties experienced with the cup product in the past
arose from the great variety of choices of $d_{\#}$, any particular choice
giving rise to artificial looking formulas''.
\end{quote}

\noindent We advise alos the reader to consult [La], [K], and [S] for details
on multiplicative rules in the intersection ring of $G_{n,k}$ in algebraic
geometry, and their history.

\bigskip

Bott-Samelson cycles provide a way to study Problem 4. To explain this we turn
back to the constructions in 4-3. We observe that

(i) The resolution $\widetilde{M}$ of $M$ are constructed from the most
familiar manifolds as

$\mathbb{R}P^{n-1}=$the real projective space of lines through the origin in
$\mathbb{R}^{n}$;

$\mathbb{C}P^{n-1}=$the real projective space of lines through the origin in
$\mathbb{C}^{n}$;

$G_{2}(\mathbb{R}^{2n})=$the Grassmannian of oriented $2$-dimensional
subspaces in $\mathbb{R}^{2n}$

\noindent and whose cohomology are well known as

\begin{center}
$H^{\ast}(\mathbb{R}P^{n-1};\mathbb{Z}_{2})=\mathbb{Z}_{2}[t]/t^{n}$;$\qquad
H^{\ast}(\mathbb{C}P^{n-1};\mathbb{Z})=\mathbb{Z}[x]/x^{n}$;

$H^{\ast}(G_{2}(\mathbb{R}^{2n});\mathbb{Z})=\{%
\begin{array}
[c]{c}%
\mathbb{Z}[y,v]/<x^{n}-2x\cdot v,v^{2}>\text{ if }n\equiv1\operatorname{mod}%
2;\ \ \ \quad\quad\quad\ \ \\
\mathbb{Z}[y,v]/<x^{n}-2x\cdot v,v^{2}-x^{n-1}\cdot v>\text{if }%
n\equiv0\operatorname{mod}2\
\end{array}
$
\end{center}

\noindent where

\begin{quote}
(a) $t$($\in H^{1}(\mathbb{R}P^{n-1};\mathbb{Z}_{2})$) is the Euler class for
the canonical real line bundle over $\mathbb{R}P^{n-1}$;

(b) $x$($\in H^{2}(\mathbb{C}P^{n-1};\mathbb{Z})$) is the Euler class of the
real reduction for the canonical complex line bundle over $\mathbb{C}P^{n-1}$;

(c) $y$($\in H^{2}(G_{2}(\mathbb{R}^{2n});\mathbb{Z})$) is the Euler class of
the canonical oriented real $2$-bundle $\gamma$ over $G_{2}(\mathbb{R}^{2n})$,
and where if $s\in H^{2n-2}(G_{2}(\mathbb{R}^{2n});\mathbb{Z})$ is the Euler
class for the orthogonal complement $\nu$ of $\gamma$ in $G_{2}(\mathbb{R}%
^{2n})\times\mathbb{R}^{2n}$, then
\end{quote}

\begin{center}
$v=\frac{1}{2}(y^{n-1}+s)\in H^{2n-2}(G_{2}(\mathbb{R}^{2n});\mathbb{Z}%
)$\footnote{The ring $H^{\ast}(G_{2}(\mathbb{R}^{2n});\mathbb{Z})$ is torsion
free. The class $y^{n-1}+s$ is divisible by $2$ because of $w_{2n-2}%
(\nu)\equiv s\equiv y^{n-1}$ $\operatorname{mod}2$, where $w_{i}$ is the
$i^{th}$ Stiefel-Whitney class.}.
\end{center}

(ii) the manifolds $\widetilde{M}$ are simpler than $M$ either in terms of
their geometric formation or of their cohomology

$H^{\ast}(\widetilde{M};\mathbb{Z})=\mathbb{Z}_{2}[t_{1},\cdots,t_{n^{\prime}%
}]/<t_{i}^{n},1\leq i\leq n^{\prime}>$ if $M=SO(n)$;

$H^{\ast}(\widetilde{M};\mathbb{Z})=\mathbb{Z}[x_{1},\cdots,x_{k}%
]/<p_{i},1\leq i\leq k>$ if $M=G_{n,k}$; and

$H^{\ast}(\widetilde{M};\mathbb{Z})=\mathbb{Q}[y_{1},\cdots,y_{n}%
]/<e_{i}(y_{1}^{2},\cdots,y_{n}^{2}),1\leq i\leq n-1;y_{1}\cdots y_{n})>$

\noindent if $\widetilde{M}=\mathbb{C}S_{n}$, where $p_{i}$ is the component
of the formal polynomial

\begin{center}
\ $\prod\limits_{1\leq s\leq i}(1+x_{s})^{-1}$
\end{center}

\noindent in degree $2(n-i+1)$ (cf. [D$_{3}$, Theorem 1]), and where
$e_{j}(y_{1}^{2},\cdots,y_{n}^{2})$ is the $j^{th}$ elementary symmetric
function in the $y_{1}^{2},\cdots,y_{n}^{2}$.

(iii) Bott-Samelson cycles on $M$ can be obtained by restricting
$h:\widetilde{M}\rightarrow M$ to appropriate subspaces of $\widetilde{M}$
(cf. Construction 2).

\bigskip

One can infer from (iii) the following result.

\textbf{Theorem 7.} \textsl{The induced ring map }$h^{\ast}:H^{\ast
}(M;\mathbb{Z}$\textsl{ or }$Z_{2})\rightarrow H^{\ast}(\widetilde
{M};\mathbb{Z}$\textsl{ or }$\mathbb{Z}_{2})$\textsl{ is injective. Furthermore}

\textsl{(1) if }$M=SO(n)$\textsl{, then}

\begin{center}
$h^{\ast}(\Omega(I))=m_{I}(t_{1},\cdots,t_{n^{\prime}})$\textsl{,}
\end{center}

\noindent\textsl{where }$m_{I}(t_{1},\cdots,t_{n^{\prime}})$\textsl{ is }the
monomial symmetric function\textsl{ in }$t_{1},\cdots,t_{n^{\prime}}$\textsl{
associated to the partition }$I$\textsl{ ([D}$_{2}$\textsl{]);}

\textsl{(2) if }$M=G_{n,k}$\textsl{, then}

\begin{center}
$h^{\ast}(\Omega(I))=S_{I}(x_{1},\cdots,x_{k})$\textsl{,}
\end{center}

\noindent\textsl{where }$S_{I}(x_{1},\cdots,x_{k})$\textsl{ is }the Schur
Symmetric function\textsl{ in }$x_{1},\cdots,x_{k}$\textsl{ associated to the
partition }$I$\textsl{ ([D}$_{1}$\textsl{]);}

\textsl{(3) if }$M=CS_{n}$\textsl{, then}

\begin{center}
$h^{\ast}(\Omega(I))=P_{I}(y_{1},\cdots,y_{n})$\textsl{,}
\end{center}

\noindent\textsl{where }$P_{I}(y_{1},\cdots,y_{n})$\textsl{ is }the Schur P
symmetric function\textsl{ in }$y_{1},\cdots,y_{n}$\textsl{ associated to the
partition }$I$\textsl{.}$\square$

(For definitions of these symmetric functions, see [Ma]).

\bigskip

Indeed, in each case concerned by Theorem 7, it can be shown that the
$\Omega(I)$ are the Schubert classes [Ch, BGG]. \

It was first pointed out by Giambelli [G$_{1}$,G$_{2}$] in 1902 (see also
Lesieur [L] or Tamvakis [T] in this volume) that multiplicative rule of
Schubert classes in $G_{n,k}$ formally coincides with that of Schur functions,
and by Pragacz in 1986 that multiplicative rule of Schubert classes in
$\mathbb{C}S_{n}$ formally agree with that of Schur P functions [P, \S6]. Many
people asked why such similarities could possibly occur [S]. For instance it
was said by C. Lenart [Le] that

\begin{quote}
``No good explanation has been found yet for the occurrence of Schur functions
in both the cohomology of Grassmanian and representation theory of symmetric groups''.
\end{quote}

\noindent Theorem 7 provides a direct linkage from Schubert classes to
symmetric functions. It is for this reason combinatorial rules for multiplying
symmetric functions of the indicated types (i.e. the \textsl{monomial
symmetric functions, Schur symmetric functions }and\textsl{ Schur P symmetric
functions}) correspond to the intersection products of Schubert varieties in
the spaces $M=SO(n)$, $G_{n,k}$ and $\mathbb{C}S_{n}$.

\textbf{Remark.} A link between representations and homogeneous spaces is
furnished by Borel [B].

\begin{center}
\textbf{4--5. A concluding remark}
\end{center}

Bott is famous for his periodicity theorem, which gives the homotopy groups of
the matrix groups $O(n;\mathbb{F})$ with $\mathbb{F=R},\mathbb{C}$ or
$\mathbb{H}$ in the stable range. However, this part of Bott's work was
improved and extended soon after its appearance [Ke], [HM], [AB].

It seems that the idea of Morse functions of Bott-Samelson type appearing
nearly half century ago [BS$_{1}$, BS$_{2}$] deserves further attention.
Recently, an analogue of Theorem 7 for the induced homomorphism

\begin{center}
$g_{w}^{\ast}:H^{\ast}(G/T)\rightarrow H^{\ast}(\Gamma_{w})$
\end{center}

\noindent of the Bott-Samelson cycle $g_{w}:\Gamma_{w}\rightarrow G/T$ (cf.
Theorem 5) is obtained in [D$_{4}$, Lemma 5.1], from which the multiplicative
rule of Schubert classes and the Steenrod operations on Schubert classes in a
generalized flag manifold $G/H$ [Ch, BGG] have been determined [D$_{4}$],
[DZ$_{1}$], [DZ$_{2}$], where $G$ is a compact connected Lie group, and where
$H\subset G$ is the centralizer of a one-parameter subgroup in $G$.

\begin{center}
\textbf{References}
\end{center}

[AB] M. Atiyah and R. Bott, On the periodicity theorem for complex vector
bundles, Acta Mathematica, 112(1964), 229-247.

[B] A. Borel,Borel, Armand Sur la cohomologie des espaces fibr\'{e} principaux
et des espaces homog\'{e}nes de groupes de Lie compacts, Ann. of Math. (2) 57,
(1953). 115--207.

[BGG] I. N. Bernstein, I. M. Gel'fand and S. I. Gel'fand, Schubert cells and
cohomology of the spaces G/P, Russian Math. Surveys 28 (1973), 1-26.

[Br] M. Brion, Lectures on the geometry of flag varieties, \textsl{this
volume}.

[BS$_{1}$] R. Bott and H. Samelson, The cohomology ring of G/T, Nat. Acad.
Sci. 41 (7) (1955), 490-492.

[BS$_{2}$] R. Bott and H. Samelson, Application of the theory of Morse to
symmetric spaces, Amer. J. Math., Vol. LXXX, no. 4 (1958), 964-1029.

[Ch] C. Chevalley, Sur les D\'{e}compositions Celluaires des Espaces G/B, in
Algebraic groups and their generalizations: Classical methods, W. Haboush ed.
Proc. Symp. in Pure Math. 56 (part 1) (1994), 1-26.

[D] J. Dieudonn\'{e}, A history of Algebraic and Differential Topology,
1900-1960, Boston, Basel, 1989.

[D$_{1}$] H. Duan, Morse functions on Grassmanian and Blow-ups of Schubert
varieties, Research report 39, Institute of Mathematics and Department of
Mathematics, Peking Univ., 1995.

[D$_{2}$] H. Duan, Morse functions on Stiefel manifolds Via Euclidean
geometry, Research report 20, Institute of Mathematics and Department of
Mathematics, Peking Univ., 1996.

[D$_{3}$] H. Duan, Some enumerative formulas on flag varieties, Communication
in Algebra, 29 (10) (2001), 4395-4419.

[D$_{4}$] H. Duan, Multiplicative rule of Schubert classes, to appear in
Invent. Math. (cf. arXiv: math. AG/ 0306227).

[DZ$_{1}$] H. Duan and Xuezhi Zhao, A program for multiplying Schubert
classes, arXiv: math.AG/0309158.

[DZ$_{2}$] H. Duan and Xuezhi Zhao, Steenrod operations on Schubert classes,
arXiv: math.AT/0306250.

[Eh] C. Ehresmann, Sur la topologie de certains espaces homog\`{e}nes, Ann. of
Math. 35(1934), 396-443.

[G$_{1}$] G. Z. Giambelli, Risoluzione del problema degli spazi secanti, Mem.
R. Accad. Sci. Torino (2)52(1902), 171-211.

[G$_{2}$] G. Z. Giambelli, Alcune propriet\`{a} delle funzioni simmetriche
caratteristiche, Atti Torino 38(1903), 823-844.

[Han] H.C. Hansen, On cycles in flag manifolds, Math. Scand. 33 (1973), 269-274.

[H] M. Hirsch, Differential Topology, GTM. No.33, Springer-Verlag, New
York-Heidelberg, 1976.

[HM] C. S. Hoo and M. Mahowald, Some homotopy groups of Stiefel manifolds,
Bull.Amer. Math. Soc., 71(1965), 661--667.

[HPT] W. Y. Hsiang, R. Palais and C. L. Terng, The topology of isoparametric
submanifolds, J. Diff. Geom., Vol. 27 (1988), 423-460.

[K] S. Kleiman, Problem 15. Rigorous fundation of the Schubert's enumerative
calculus, Proceedings of Symposia in Pure Math., 28 (1976), 445-482.

[Ke] M.A. Kervaire, Some nonstable homotopy groups of Lie groups, Illinois J.
Math. 4(1960), 161-169.

[La] A. Lascoux, Polyn\^{o}mes sym\'{e}triques et coefficients d'intersection
de cycles de Schubert. C. R. Acad. Sci. Paris S\'{e}r. A 279 (1974), 201--204.

[L] L. Lesieur, Les problemes d'intersections sur une variete de Grassmann, C.
R. Acad. Sci. Paris, 225 (1947), 916-917.

[Le] C. Lenart, The combinatorics of Steenrod operations on the cohomology of
Grassmannians, Advances in Math. 136(1998), 251-283.

[Ma] I. G. Macdonald, Symmetric functions and Hall polynomials, Oxford
Mathematical Monographs, Oxford University Press, Oxford, second ed., 1995.

[M] C. Miller, The topology of rotation groups, Ann. of Math., 57(1953), 95-110.

[M$_{1}$] J. Milnor, Lectures on the h-cobordism theorem, Princeton University
Press, 1965.

[M$_{2}$] J. Milnor, Morse Theory, Princeton University Press, 1963.

[M$_{3}$] J. Milnor, Differentiable structures on spheres, Amer. J. Math.,
81(1959), 962-972.

[MS] J. Milnor and J. Stasheff, Characteristic classes, Ann. of Math. Studies
76, Princeton Univ. Press, 1975.

[P] P. Pragacz, Algebro-geometric applications of Schur S- and Q-polynomi-als,
Topics in invariant Theory (M.-P. Malliavin, ed.), Lecture Notes in Math.,
Vol. 1478, Spring-Verlag, Berlin and New York, 1991, 130-191.

[Sch] H. Schubert, Kalk\"{u}l der abz\"{a}hlende Geometric, Teubner, Leipzig, 1879.

[S] R. P. Stanley, Some combinatorial aspects of Schubert calculus, Springer
Lecture Notes in Math. 1353 (1977), 217-251.

[St] N. E. Steenrod and D. B. A. Epstein, Cohomology Operations, Ann. of Math.
Stud., Princeton Univ. Press, Princeton, NJ, 1962.

[T] H. Tamvakis, Gromov-Witten invariants and quantum cohomology of
Grassmannians, \textsl{this volume}.

[VD] A.P. Veselov and I.A. Dynnikov, Integrable Gradient flows and Morse
Theory, Algebra i Analiz, Vol. 8, no 3.(1996), 78-103; Translated in St.
Petersburgh Math. J., Vol. 8, no 3.(1997), 429-446.

[Wh] J. H. C. Whitehead, On the groups $\pi_{r}(V_{n,m})$, Proc. London Math.
Soc., 48(1944), 243-291.
\end{document}